\RequirePackage{fix-cm}
\documentclass{svjour3}                     
\smartqed  
\usepackage[final]{graphicx}
\usepackage[utf8]{inputenc}
\usepackage[english]{babel}
\usepackage[T1]{fontenc}
\usepackage[T1]{fontenc}
\usepackage{lmodern}
\usepackage{graphicx}
\usepackage{stmaryrd} 
\usepackage{epstopdf}
\usepackage{amsfonts}
\usepackage{dsfont}
\usepackage{amsmath}
\usepackage{array}
\usepackage{amssymb}
\usepackage{enumerate}
\usepackage{euscript}
\usepackage{ulem} 
\usepackage{manfnt}
\usepackage[english]{babel}

\newcommand{\ppe}{\vspace*{0.1cm}}
\newcommand{\pe}{\vspace*{0.2cm}}

\newcommand{\bea}{\begin{array}}
\newcommand{\ena}{\end{array}}

\def\[ent{[\hskip -1.5pt [}
\def\]ent{]\hskip -1.5pt ]}

\newcommand{\ie}[2]{\,\llbracket #1,\,#2 \rrbracket} 

\newcommand \R {\mathbb{R}}
\newcommand \C {\mathbb{C}}
\newcommand \N {\mathbb{N}}
\newcommand \Z {\mathbb{Z}}

\renewcommand{\P}{\mathbb{P}}
\newcommand \E {\mathbb{E}}

\usepackage{tikz, ifthen}
	\usetikzlibrary{arrows,shapes,positioning,automata,decorations.markings}
	\usepackage{float}
\tikzstyle directed=[postaction={decorate,decoration={markings,
    mark=at position .65 with {\arrow{latex}}}}]
\usepackage[margin=20pt,format=plain,small,justification=centering]{caption}

\usepackage{mathtools}

\usepackage{theorem}
\theorembodyfont{\upshape}
\newtheorem{asn}{Assumption}
\newtheorem{lea}{Lemma}
\newtheorem{prn}{Proposition}
\newtheorem{thth}{Theorem}
\newtheorem{coe}{Corollary}
\begin{document}

\title{A stochastic SIR model on a graph with epidemiological and population dynamics occurring over the same time scale
}
%



\author{Pierre Montagnon}


\institute{Pierre Montagnon \at
              CMAP, École Polytechnique, Route de Saclay, 91128 Palaiseau Cedex, France. \\
              MaIAGE, INRA, Université Paris-Saclay, 78350 Jouy-en-Josas, France. \\
              \email{pierre.montagnon@polytechnique.edu}           
}

\date{Received: date / Accepted: date}

\maketitle

\begin{abstract}
We define and study an open stochastic SIR (Susceptible -- Infected -- Removed) model on a graph in order to describe the spread of an epidemic on a cattle trade network with epidemiological and demographic dynamics occurring over the same time scale. Population transition intensities are assumed to be density-dependent with a constant component, the amplitude of which determines the overall scale of the population process. Standard branching approximation results for the epidemic process are first given, along with a numerical computation method for the probability of a major epidemic outbreak. This procedure is illustrated using real data on trade-related cattle movements from a densely populated livestock farming region in western France (Finistère) and epidemiological parameters corresponding to an infectious epizootic disease. Then we exhibit an exponential lower bound for the extinction time and the total size of the epidemic in the stable endemic case as a scaling parameter goes to infinity using results inspired by the Freidlin-Wentzell theory of large deviations from a dynamical system. 
\end{abstract}

\section{Introduction}
\label{intro}\sectionmark{Introduction}

\indent \indent Animal movements are a major vector of epidemic propagation between cattle holdings at large spatial scales. The large amount of data collected by European authorities over the last two decades \cite{Ver,Dut} makes it possible to track the position over time of every single piece of cattle within national territories, allowing to design and calibrate models for pathogen spread \cite{Hos,Perra}. For some diseases, modelling the propagation of an epidemic on a cattle trade network requires taking into account demographic and epidemiological dynamics occurring at the same time scale. Mathematically, this involves coupling epidemiological multitype stochastic processes (Chapter 6 of \cite{AB}, \cite{BC93,Cla96}) with demographic models \cite{ONe,Nas,AB00} with births, movements between nodes and deaths (that we should refer to as \textit{open} demographic models). In the SIR (Susceptible -- Infected -- Removed) case, it has been shown \cite{Nas,AB00,VH95,VH97} that such models differ from their demography-free counterparts in that they allow for endemicity, that is, for the persistence of an epidemic over a given threshold for a long period of time. The probability of a major epidemic outbreak (Chapter 4 of \cite{AB}, \cite{Cla96,Nea12}) to occur, the extinction time and the total size of an epidemic (that is, the total number of individuals infected during the course of the epidemic) in the case of a major outbreak are therefore essential quantities of interest in the study of the epidemic process, especially from a control perspective.\pe

\indent \indent There is a great amount of literature about dynamical epidemics on a contact network. While many authors consider individuals as nodes (e.g. \cite{BN,DF15}), some think of nodes --- or \textit{types} --- as subpopulations \cite{BC93,BC95,Nea12}, which results in coupling epidemiological and metapopulation models. In the latter case, intra-nodal population dynamics related to births or deaths are seldom taken into account, thus yielding fixed-size models in which individuals may \cite{Nea12,BC93} or may not \cite{BC95} move across nodes. According to cases, infectives may or may not make infectious contacts with individuals from other nodes.

Some authors consider open demographic dynamics coupled with epidemiological processes but mostly deal with single-type models \cite{Nas,AB00,ONe} --- that is, with one single, uniformly mixing population. Moreover, although most of such models either rely on a population process with density-proportional \cite{ONe,CG} or constant \cite{Nas,AB00,VH95,VH97} entry rates, the specific modelling of livestock demographic dynamics requires to account of both Malthusian population growth and immigration into the area under study, so we want to introduce both a population-proportional birth component and a constant immigration component for entry rates. The specification we choose is reminiscent of a few previous papers. \cite{BS} performed a numerical analysis of a single-type model ($n=1$) only differing from ours by the possibility of vertical health status transmission. Close deterministic counterparts for our stochastic model can be found in \cite{McK,LS,TL,MEK,NGT}.\pe 

\indent \indent In the present paper, we  model cattle holdings as the nodes of a directed graph, the edges of which are formed by trade paths between nodes. We consider an open, multitype population process on this graph with state-dependent, affine birth rates. Such a process, that is reminiscent of stochastic metapopulation models \cite{Verb}, is a multitype continuous-time branching process \cite{Mode,Ath68,AthNey} with immigration, which leads to strong stability properties under a mild subcriticality assumption. Our second modelling step is to define an individual-based SIR multitype process \cite{BC93,BC95,Nea12} within the network nodes. Each node will be endowed with its proper epidemic-related parameters and may receive or send susceptible, infective or removed individuals. We consider only one local level of mixing: infectives may only make infectious contacts with individuals from the same node. Movements between nodes are therefore the only reason for the spread of the epidemic across nodes (this setting being referred to as a \textit{dynamical epidemic model} in \cite{Nea12}).\pe

Our next step (Section \ref{sec_major_outbreak_th}) is to generalize closed-population branching approximation results \cite{AB,Cla96,Nea12} to our setting as a population scaling parameter goes to infinity. We define the basic reproduction number $R_0$ for the epidemic process and compute the associated major outbreak probability, using and refining in our Markovian framework a procedure put forward by \cite{Nea12}. This allows for a numerical application using data on the Finistère cattle trade network from the French National Identification Database (BDNI) in Section \ref{sec_major_outbreak_num}.

\pe 

Finally, we will discuss the behavior of the epidemic process in the case of a major outbreak (Section 4). We will exhibit a lower bound for the maximal number of infectives during the course of the epidemic and use it to derive an exponential lower bound for the extinction time and the total size of the epidemic in the case where the associated dynamical system admits an endemic equilibrium using an adaptation of the Freidlin-Wentzell theory for large deviations to Poisson perturbations. The latter bound is the main result of this paper, and we have good hints to believe that it gives the right magnitude order for the epidemic extinction time and total size. It differs fundamentally from those derived for fixed-size SIR models (where the total size of the epidemic has to be lower than the population total size and its extinction time is proportional to a logarithm of the scaling parameter, see Chapter 4 of \cite{AB}) and illustrates the role of population renewal in the persistence of an endemic disease.

\pe 

Let $n\geqslant 1$. In the rest of the paper, we shall write $\cdot$ for the usual scalar product, $\Vert \cdot \Vert_i$ for $l_1$-norms if $i=1$, Euclidean norms if $i=2$ and uniform norms if $i=\infty$, $\mathcal{B}_i(x,\delta)$ for the open ball with center $x$ and radius $\delta$ for $\Vert\cdot \Vert_i$ and $\text{diag}(\lambda_1,\ldots,\lambda_n)$ for the diagonal matrix of $\mathcal{M}_n(\R)$ with diagonal coefficients $\lambda_1,\ldots,\lambda_n$. Finally, the $i$-th coordinate of vector $x$ will be denoted as $x_i$. All processes will be assumed to be defined on a unique measurable space $(\Omega,\mathcal{A})$.

\section{Model setting}
\label{sec_data_model}

\begin{figure}[!h]
\begin{center}
\includegraphics[scale=.28]{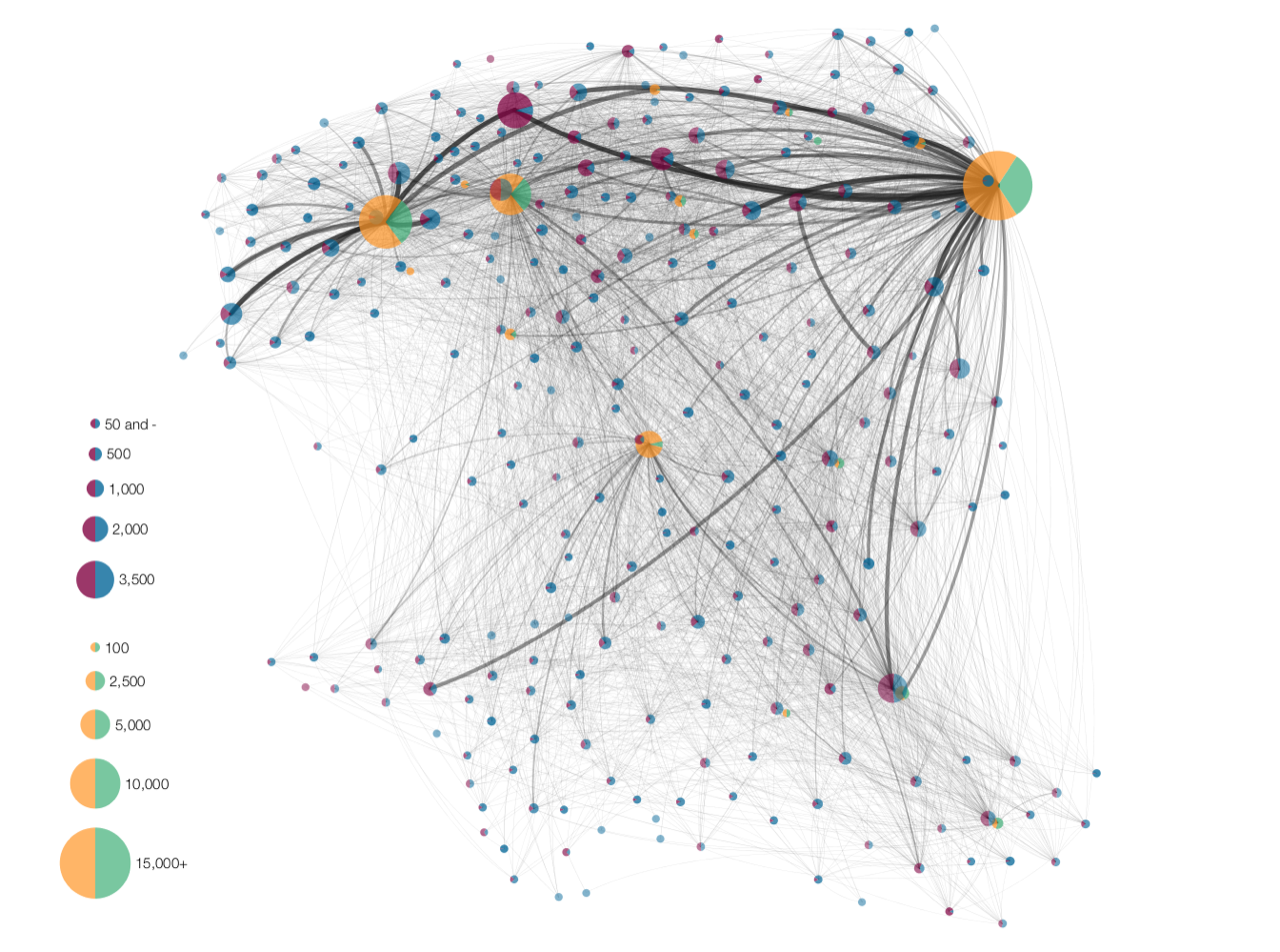}
\caption{\label{fig_finistère}Geographical location of exchanging holdings of Finistère (a Western region of France densely populated with cattle) in 2015, aggregated at the commune level. The size of a node represents its mean population over the year. Colours on the nodes represent the share of buying (orange for operators and red for farms) and selling (green for operators and blue for farms). The edges represent the existence of a movement between two holdings in the dataset, with their width indicative of the observed trading volume along this edge. The picture is courtesy of Gaël Beaunée.}
\end{center}
\end{figure}

\indent \indent Modelling choices will be guided by empirical observations from the BDNI. This database traces back the path of every piece of cattle within the French territory \cite{Dut,Rau} between 2005 and 2016. We extracted information corresponding to a densely populated livestock farming region in western France (Finistère), for 2015. The 2015 Finistère exchange subnetwork we consider in this paper (see Figure \ref{fig_finistère}) contains 4,183 vertices and 10,036 edges. It consists of 4,163 farms (with internal population dynamics) and commercial operators, 3 markets and 17 assembly centers (without internal population dynamics) exchanging a total of 118,311 animals, receiving 55,325 from the outside and sending 241,747 to nodes outside of Finistère or slaughterhouses (these movements being considered as deaths). The average total population on the network over the year is 424,385. These numbers are quite large, which motivates the introduction of some scaling parameter.\pe 

\subsection{The population model}
\label{subsec_pop_mod}

\indent \indent We first define a $\Z_+^n$-valued population process $(X^N_t)_{t\geqslant 0}$ modelling the population dynamics ($n$ being the number of nodes in the graph and $N$ a scaling parameter) and state a condition (Assumption \ref{assump_sub}) ensuring its convergence in distribution, as $t\to \infty$, to an invariant probability measure at a geometric speed (Proposition \ref{th_stabpop}). We then investigate the time for the scaled population process to deviate by any given fixed distance from the corresponding deterministic model, described by the solution of a linear ODE (Propositions \ref{lgn_pop} and \ref{thm_fw_pop}).\pe

\indent We model the population dynamics using a $\Z_+^n$-valued multitype continuous-time branching process (BP) with immigration. At any population state $x=(x_1,\ldots,x_n)\in \Z_+^n$, the inflow rate in node $i$ is $NB_i+b_i x_i$ and the death rate in node $i$ is $d_i x_i$ with $b_i,B_i,d_i\in \R_+$. The $b_ix_i$ might be considered as Malthusian birth rates and the $NB_i$ as constant immigration rates. 

Data suggest that the temporal rate of transfers between agents also is size-dependent. Now their amplitude is bounded because of transportation constraints and does not vary much empirically, so we set it to be unitary and define the transfer rate from node $i$ to node $j$ at population state $x$ as $\theta_{i,j}x_i$ with $\theta_{i,j}\geqslant 0$. 

\pe

Let $x(0)\in \R_+^n$. For any $N>0$, define $(X^N(t))_{t\geqslant 0}$ as a $\Z_+^n$-valued jump process with initial value $X^N(0)=\lfloor N x(0) \rfloor$ and the following transition rates under $\P$:
\begin{equation}\label{tableau_taux_trans}\begin{array}{cc}
\text{Transition} & \text{Rate at state }x \\
x\to x+e_i & NB_i+b_ix_i \\
x\to x-e_i & d_ix_i \\
x\to x-e_i+e_j & \theta_{i,j}x_i
\end{array}\end{equation}
where $(e_1,\ldots,e_n)$ is the canonical basis of $\R^n$. All $X^N$ are assumed to be built using a single set of independent, homogeneous Poisson processes with rate $1$ and random time changes (see page 326 of \cite{EK} and the proof of Theorem \ref{thm_approx_sir_bp} below). If $x\in \Z_+^n$ and if $\mu$ is a probability distribution on $\Z_+^n$, we will use the classical notations $\P_{x}$ and $\P_{\mu}$ to denote probabilities on $(\Omega,\mathcal{A})$ under which, for any $N>0$, $X^N(0)$ has respective distributions $\delta_x$ and $\mu$. Associated expectations will be denoted by $\E_x$ and $\E_{\mu}$.

\indent It can easily be shown that $X^N$ is non-explosive, that is, $X^N(t)$ is finite for all $t\in \R_+$ with probability $1$ for any initial value $x(0)\in \Z_+^n$.

\pe

\indent From now on, we will assume that the directed graph with vertex set $\{1,\ldots,n\}$ and edge set $\{(i,j)\mid \theta_{i,j}>0\}$ is fully connected, so that for any $i,j$, any individual born in node $i$ can get to node $j$ with positive probability during its lifetime.

\pe 

\indent In our framework, immigration in node $i$ occurs at rate $NB_i$, and individuals in node $i$ give birth at rate $b_i$, die at rate $d_i$ and move to node $j$ at rate $\theta_{i,j}$, independently from other individuals in the network. The lineage of a single individual is a branching process with transition rates given by (\ref{tableau_taux_trans}) with the $B_i$ replaced by $0$. Once they enter the system, individuals give birth to independent lineages that do not interact, so the number of such lineages and therefore the typical size of the population process is proportional to the mean vector of immigrants per unit of time (that is, $B=(B_1,\ldots,B_n)$).\pe 

We wish to model populations that are stable over time, so it is clear that this branching process must not tend to infinity with positive probability, that is, we do not want it to be supercritical (see however \cite{BT} for the study of an epidemic process within a growing population modelled using a supercritical branching process). We impose a slightly stronger condition --- subcriticality --- to make sure that the first moment of the population process does not go to infinity. From now on, we will assume the following condition to hold.

\begin{asn}[Subcriticality of the immigration-free population BP]\label{assump_sub}
The eigenvalues of
$$A=\begin{pmatrix} b_1-d_1-\sum_{j\neq 1}\theta_{1,j} & \theta_{2,1} & \cdots & \theta_{n,1}\pe \\ \theta_{1,2} & b_2-d_2-\sum_{j\neq 2}\theta_{2,j} & \ddots & \vdots\pe  \\ \vdots & \ddots & \ddots & \theta_{n,n-1}\pe  \\ \theta_{1,n} & \cdots & \theta_{n-1,1} & b_n-d_n-\sum_{j\neq n}\theta_{n,j}\end{pmatrix}$$
have negative real parts.
\end{asn}

 We may then state our first convergence result for the population process. Although it is a continuous-time version of standard results for discrete-time multitype branching processes with immigration, we could not find the exact same statement in the existing literature. We present its proof in the Appendix.

\begin{prn}\label{th_stabpop}
\textit{Let $N>0$. $X^N$ is positive recurrent and $(\Vert \cdot\Vert_1+1)-$exponentially ergodic, that is, the invariant probability $\pi$ of $X^N$ has a finite first-order moment and there exist $\beta\in (0,1)$ and $\lambda\in (0,+\infty)$ independent from the choice of $x(0)$ such that:
$$\forall t\geqslant 0,\quad \sup_{|g|\leqslant \Vert \cdot \Vert_1+1}\left|\E_{x(0)}(g(X^N(t)))-\int g \mathrm{d}\pi\right|\leqslant \lambda(\Vert x(0) \Vert_1+1)\beta^t,$$
where the supremum is taken over all measurable functions $g:\Z_+^n\to \C$ such that $|g(x)|\leqslant \Vert x \Vert_1+1$ for all $x\in \Z_+^n$. Moreover,
$$\lim_{t\to +\infty}\E_{x(0)}(X^N(t))=\int x\mathrm{d}\pi(x)=-NA^{-1}B.$$}
\end{prn}

Just as expected, the limiting average population size $-NA^{-1}B$ is proportional to the scaling factor $N$.

\begin{figure}[!h]
\hspace{-3mm}\begin{center}
\includegraphics[scale=.62]{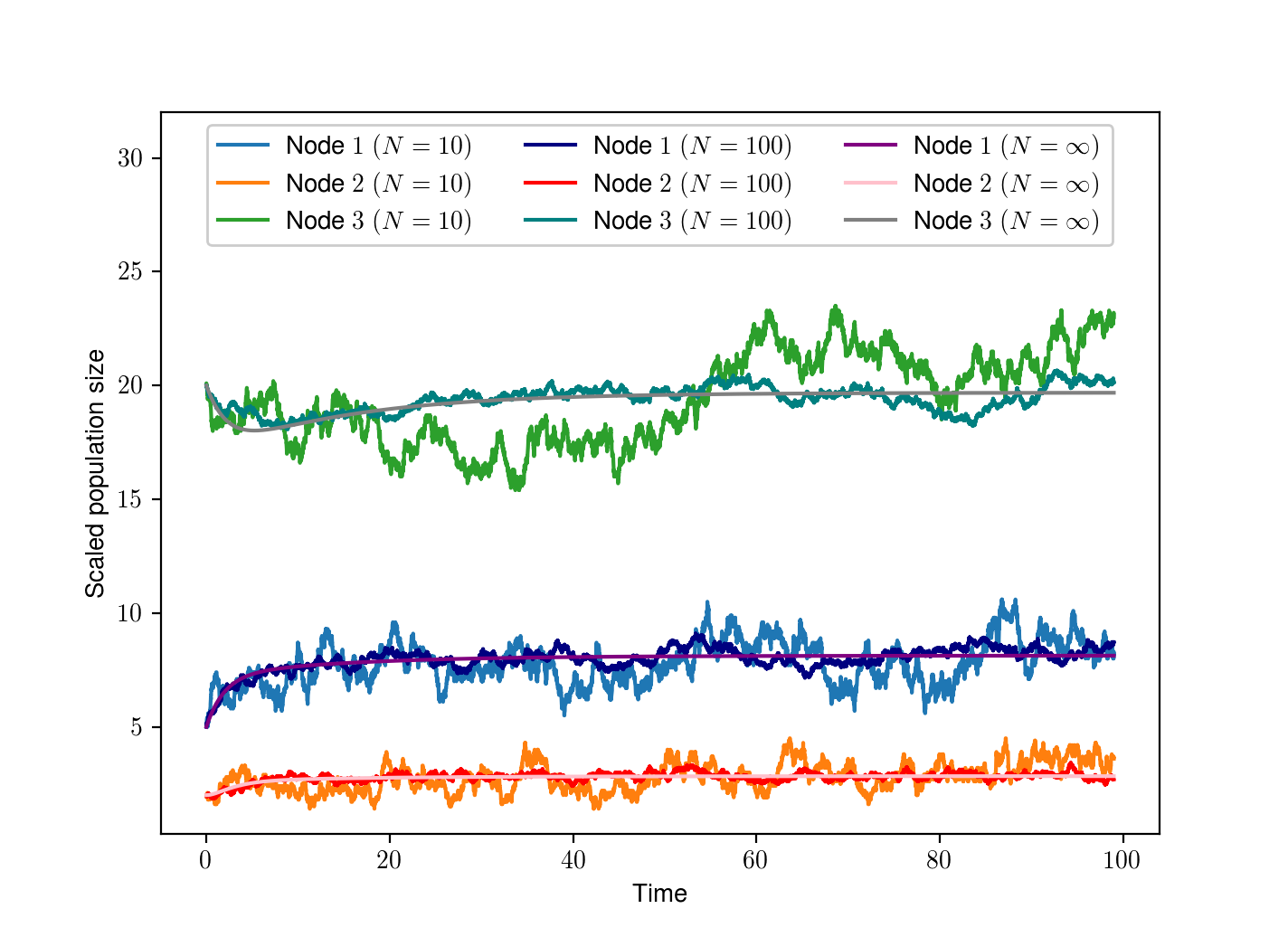}
\end{center}
\caption{\label{fig_smooth_pop}Simulated values of $X^N_t/N$ ($t\in [0,100]$) for $N=10$ and $N=100$ and limiting deterministic process with arbitrary parameter values $n=3$, $\tilde{x}_0=(5,2,20)$, $\tilde{B}=(1,.5,.5)$, $d=(.1,.2,.1)$, $b=(.1,.05,.02)$, $\theta_{1,2}=\theta_{1,3}=.2$, $\theta_{2,1}=.1$, $\theta_{2,3}=.5$, $\theta_{3,1}=.1$ and $\theta_{3,2}=0$. The time unit is arbitrary.}
\end{figure}

We are interested in the behavior of our process as $N$ tends to infinity. 
The mean equilibrium value $z^*=-A^{-1}B$ of $X^N/N$ does not depend on $N$, and as $N$ grows we expect the scaled superimposition of independent lineages to get smoother (see Figure \ref{fig_smooth_pop}, where the parameter values are chosen arbitrarily to make the figure easy to read and assimilate). Standard results on scaling limits of density dependent population processes can be found in \cite{EK}, Chapter 11, and show that we may indeed approximate $X^N/N$ using Brownian deviations from a deterministic process on finite time intervals as $N$ tends to infinity. The following law of large numbers is a direct consequence of Theorem 2.1 from Chapter 11 of \cite{EK}.

\begin{prn}\label{lgn_pop}
\textit{Define $z$ as the solution of the Cauchy problem $z'=Az+B$ with $z(0)=x(0)$, that is,
$$\begin{array}{cccc}z:&\R&\longrightarrow &\R^n \\ &t&\longmapsto&e^{tA}(A^{-1}B+x(0))-A^{-1}B\end{array}~.$$
For any $T\geqslant 0$,
$$\P\left(\lim_{N\to +\infty}\sup_{t\in [0,T]}\left\Vert \frac{X^N(t)}{N}-z(t)\right\Vert_{\infty}=0\right)=1.$$}
\end{prn}

The deterministic process $z=(z(t))_{t\geqslant 0}$ quickly converges to its equilibrium value $z^*:=-A^{-1}B$. When coupling further dynamics with the population process, it is therefore common to consider that the latter starts from some point close to $z^*$. The following result provides bounds for population fluctuations over very large time intervals for such initial conditions, which will prove useful to describe the early phase of the epidemic in the next section. Its proof is based on a Freidlin-Wentzell-type results on large deviations from a deterministic approximation given in \cite{ParBrice} (see also \cite{KP,PB}) and can be found in the Appendix.

\begin{prn}\label{thm_fw_pop}
\textit{\textit{Let $\varepsilon\in (0,\Vert z^*\Vert_{\infty})$ and denote by $\tau^N_{\varepsilon}$ the exit time of the ball $\mathcal{B}_{\infty}(z^*,\varepsilon)$ by $X^N/N$. There exists $\alpha_0>0$ such that for any $\alpha>0$:}
$$\forall x\in \mathcal{B}_{\infty}(z^*,\varepsilon),\quad \lim_{N\to +\infty}\P_{\lfloor Nx \rfloor}\left(e^{(\alpha_0-\alpha)N}<\tau^N_{\varepsilon}<e^{(\alpha_0+\alpha)N}\right)=1.$$}
\end{prn}

The constant $\alpha_0$ above is the \textit{exit cost} from $\mathcal{B}_2(z^*,\varepsilon)$ starting from $z^*$ for the dynamical system $y'=Ay+B$ and the Poisson perturbation considered, that is, the minimal value of the quasipotential for this system and the perturbation with respect to $z^*$ on the boundary of this ball (see Chapter V of \cite{FW}, see \cite{ParBrice}, and see the proof of Proposition \ref{thm_fw_pop} for an expression of the quasipotential).


\subsection{The epidemic model}\label{sec_epid}

\indent \indent \indent We now define a stochastic SIR model for the spread of an epidemic within and between nodes, then generalize a standard finite-time convergence result to a branching process as $N$ goes to infinity.\pe

\indent \indent For any $i$ and $t$, the population $X^N_i(t)$ of the $i$-th node at time $t$ is divided into three subpopulations: $S^N_i(t)$ susceptibles, $I^N_i(t)$ infectives and $R^N_i(t)$ removed individuals. Movements between nodes, births and deaths --- that is, population dynamics --- are assumed to be independent from health status. All individuals entering the system, either by birth or immigration, are assigned the susceptible status  --- so we exclude vertical disease transmission or immigration of infective individuals, although this assumption can be relaxed without technical complications. When alive and in node $i$, each infected individual makes infectious contacts, with individuals chosen independently and uniformly from node $i$, at rate $\beta_i$. Such contacts may only occur within a given node, so an infective from node $i$ cannot make an infectious contact with an individual from node $j$ if $i\neq j$. Contacted individuals that are still susceptible get infected; otherwise their status does not change. Independently, infectious individuals alive in node $i$ recover at rate $\gamma_i$ and are then removed. 
All transitions are again assumed to be realized using a single set of independent Poisson processes that is the same for all values of $N$ and random time changes (see again \cite{EK}, Chapter 11, and also \cite{Cla96} for an explicit construction).
\pe

\indent Let us write the canonical basis of $\R^{3n}$ as $(e^s_1,\ldots,e^s_n,e^i_1,\ldots,e^i_n,e^r_1,\ldots,e^r_n)$. We consider for each $N>0$ a $(\Z_+^n)^3$-valued pure jump process $(S^N(t),I^N(t),R^N(t))_{t\geqslant 0}$ defined on $(\Omega,\mathcal{A})$ and described by the following transition rates under $\P$:
\begin{equation}\label{eq_taux_poisson}\begin{array}{cc}
\text{Transition} & \text{Rate at state }x \\
(s,i,r)\to (s,i,r)+e^s_j & NB_j+b_j(s_j+i_j+r_j) \\
(s,i,r)\to (s,i,r)-e^s_j & d_js_j \\
(s,i,r)\to (s,i,r)-e^i_j & d_ji_j \\
(s,i,r)\to (s,i,r)-e^r_j & d_jr_j \\
(s,i,r)\to (s,i,r)+e^s_k-e^s_j & \theta_{j,k}s_j \\
(s,i,r)\to (s,i,r)+e^i_k-e^i_j & \theta_{j,k}i_j \\
(s,i,r)\to (s,i,r)+e^r_k-e^r_j & \theta_{j,k}r_j \\
(s,i,r)\to (s,i,r)+e^i_j-e^s_j & \beta_j \frac{i_js_j}{i_j+s_j+r_j} \\
(s,i,r)\to (s,i,r)+e^r_j-e^i_j & \gamma_j i_j \\
\end{array}\end{equation}
and such that $S^N(t)+I^N(t)+R^N(t)=X^N(t)$ for all $t\geqslant 0$.\pe


\indent Given the application we consider, we are interested in describing the behavior of an epidemic started by a small number of individuals, so we do not assume $I^N(0)$ to be proportional to $N$. We will instead consider that $\P$-almost surely
$$I^N(0)=I(0), \quad S^N(0)=X^N(0)-I(0)=\lfloor Nx(0)\rfloor-I(0) \quad \text{and} \quad R^N(0)=0$$ 
with a fixed $I(0)\in \Z_+^n\setminus\{0\}^n$.\pe

\indent It is not difficult to see that $I^N$ almost surely reaches $0$ within finite time since $X^N=S^N+I^N+R^N$ is positive recurrent and all states of $\Z_+^{3n}$ lead to $\Z^n_+\times\{0\}^n\times \Z^n_+$ with positive probability. However, we may picture situations in which the epidemic dies within its first stages and others where it spreads widely (possibly within one single node at first, then across a large portion of the graph) before going extinct, which we will refer to as epidemic major outbreaks. Our concern is to evaluate the probability for major outbreaks to happen and to quantify the extinction time and the total size\footnote{It is sufficient for our purpose to define the \textit{total size} of $I^N$ as the total number of upward jumps of $\Vert I^N\Vert_1$.} of the epidemic in such cases.\pe

\indent Theorem \ref{thm_approx_sir_bp} below is a generalization of Theorem 2.1 of \cite{Cla96} and Theorem 3.3 of \cite{ONe} to our multitype open setting with density-dependent population inflow. It relies on the idea that if $N$ is large, the early stages of the epidemic look like a branching process because the probability for the first infectives of drawing non-susceptible individuals to make contact with vanishes (see \cite{Bar55} p.141 or \cite{Ken56,Ball83}). The "birth" of an individual in this branching process corresponds to the infection of a susceptible by an infective, while the "death" of an individual means either the actual death or the recovery of the corresponding infective for the epidemic process.

\begin{thth}[Finite-time convergence to a branching process]\label{thm_approx_sir_bp}
\textit{Assuming that $(\Omega,\mathcal{A})$ is large enough, there exist a $(S^N(t),I^N(t),R^N(t))_{t\geqslant 0}$ infection jump process on $\Z_+^{3n}$ with rates given by (\ref{eq_taux_poisson}) and a $\Z_+^n$-valued continuous-time multitype branching process $(I'(t))_{t\geqslant 0}$ on $(\Omega,\mathcal{A},\P)$ with $I'(0)=I(0)$ describing a population of individuals moving from node $i$ to node $j$ at rate $\theta_{i,j}$ and, while in node $i$, giving birth at rate $\beta_i$ and dying at rate $d_i+\gamma_i$, such that for any time $T\geqslant 0$, $\P$-almost surely:
$$\exists N_0\in \Z_+^* ~:~ \forall N\geqslant N_0,\forall u\in [0,T],\quad I^N(u)=I'(u).$$}
\textit{In particular, if $\tau^N$ (resp. $Z^N$) denotes the extinction time (resp. total size) of the epidemic and $\tau'$ (resp. $Z'$) that of the branching process, then
\begin{equation}\label{conv_result}\tau^N\underset{N\to +\infty}{\longrightarrow} \tau' \quad \text{and}\quad Z^N\underset{N\to +\infty}{\longrightarrow} Z'\end{equation}
$\P$-almost surely.}
\end{thth}

\indent \indent Note that the $N\to +\infty$ approximation consists in considering that susceptibles are infinitely numerous, which explains why the birth parameters $b_i$ and $B_i$ do not enter the definition of the limiting branching process.\pe

The difference between our setting and that of \cite{Cla96} is twofold: on the one hand, we consider varying infection and recovery rates $\beta_i$ and $\gamma_i$, and on the other hand we take demographics into account by allowing for immigration, births and deaths. Yet the proof of Theorem 2.1 in \cite{Cla96} adapts well to varying infection and recovery rates, and our Proposition \ref{lgn_pop} suggests that the influence of demographics on the transition rates is negligible on finite time intervals as $N$ goes to infinity. The latter observation guides the following proof of Theorem \ref{thm_approx_sir_bp}.\pe

\textsc{Proof of Theorem \ref{thm_approx_sir_bp}} ---  For all $i$ and all $j\neq i$, let $Q^B_i$, $Q^b_i$, $Q_i^{d,1}$, $Q_i^{d,2}$, $Q_i^{d,3}$, $Q^{\beta}_i$, $Q^{\gamma}_i$, $Q^1_{i,j}$, $Q^2_{i,j}$ and $Q^3_{i,j}$ be independent unit-rate Poisson processes on $(\Omega,\mathcal{A},\P)$, and let $(A_k)_{k\geqslant 0}$ be a sequence of variables on $(\Omega,\mathcal{A},\P)$ uniformly distributed on $[0,1]$, independent from each other and from the processes introduced earlier. We define the multitype branching process $I'$ as the solution of the integral equation
\begin{multline*}I'_i(t)=I_i(0)+Q^{\beta}_i\left(\int_{0}^t \beta I'_i(u)\mathrm{d}u\right)-Q^{d,2}_i\left(\int_{0}^t d_i I'_i(u)\mathrm{d}u\right)-Q^{\gamma}_i\left(\int_{0}^t \gamma_i I'_i(u)\mathrm{d}u \right)\\ +\sum_{j\neq i}\left(Q_{j,i}^2\left(\int_0^t \theta_{j,i}I'_j(u)\mathrm{d}u\right)-Q_{i,j}^2\left(\int_0^t \theta_{i,j}I'_i(u)\mathrm{d}u\right)\right)\end{multline*}
for all $i$ and all $t$. For any $N$, we finally define $(S^N_i,I^N_i,R^N_i)$ as the solution of
\begin{multline*}S^N_i(t)=\lfloor Nx_i(0)\rfloor-I_i(0)+Q^B_i\left(tNB_i\right)+Q^b_i\left(\int_{0}^t b_i X^N_i(u)\mathrm{d}u\right)\\
-Q^{d,1}_i\left(\int_{0}^t d_i S^N_i(u)\mathrm{d}u\right)-Q^{\beta}_i\left(\int_{0}^t \beta I^N_i(u)\mathrm{1}_{\frac{S^N_i(u)^-}{X^N_i(u)}\geqslant A_{Q_i\left(\int_{0}^u \beta I'_i(s)\mathrm{d}s\right)}}\mathrm{d}u\right)\\
+\sum_{j\neq i}\left(Q^1_{j,i}\left(\int_0^t \theta_{j,i}X^N_j(u)\mathrm{d}u\right)-Q^1_{i,j}\left(\int_0^t \theta_{i,j}X^N_i(u)\mathrm{d}u\right)\right),
\end{multline*}

\begin{multline*}I^N_i(t)=I_i(0)+Q^{\beta}_i\left(\int_{0}^t \beta I^N_i(u)\mathrm{1}_{\frac{S^N_i(u)^-}{X^N_i(u)}\geqslant A_{Q_i\left(\int_{0}^u \beta I'_i(s)\mathrm{d}s\right)}}\mathrm{d}u\right)\\
-Q^{d,2}_i\left(\int_{0}^t d_i I^N_i(u)\mathrm{d}u\right)-Q^{\gamma}_i\left(\int_{0}^t \gamma_i I^N_i(u)\mathrm{d}u \right)\\
+\sum_{j\neq i}\left(Q_{j,i}^2\left(\int_0^t \theta_{j,i}I^N_j(u)\mathrm{d}u\right)-Q_{i,j}^2\left(\int_0^t \theta_{i,j}I^N_i(u)\mathrm{d}u\right)\right),\end{multline*}
and 
\begin{multline*}R^N_i(t)=Q^{\gamma}_i\left(\int_{0}^t \gamma_i I^N_i(u)\mathrm{d}u \right)-Q^{d,3}_i\left(\int_{0}^t d_i R^N_i(u)\mathrm{d}u\right)\\
+\sum_{j\neq i}\left(Q_{j,i}^3\left(\int_0^t \theta_{j,i}R^N_j(u)\mathrm{d}u\right)-Q_{i,j}^3\left(\int_0^t \theta_{i,j}R^N_i(u)\mathrm{d}u\right)\right)\end{multline*}
for all $i$ and all $t$, where $X^N:=S^N+I^N+R^N$ is a $\Z_+^n$-valued jump process with transition rates given by (\ref{tableau_taux_trans}).

It is not difficult to see that $I'$ is a $\Z_+^n$-valued branching process with the expected transition rates and that the $(S^N,I^N,R^N)$ are Poisson processes on $\Z_+^{3n}$ with transition rates given by (\ref{eq_taux_poisson}). Processes $I^N$ and $I'$ coincide up to the smallest $u\geqslant 0$ such that
$$\frac{S^N_i(u)^-}{X^N_i(u)}< A_{Q_i\left(\int_{0}^u \beta I'_i(s)\mathrm{d}s\right)}$$
for some $i$. Now if $I^+_i(u)$ denotes the total number of births in node $i$ for the branching process $I'$ up to time $u$, is it straightforward that for any fixed $u$,
$$\frac{S^N_i(u)^-}{X^N_i(u)}\geqslant \frac{X^N_i(u)-I_i(0)-I_i^+(u)}{X^N_i(u)}\underset{N\to +\infty}{\longrightarrow} 1,$$
so for any fixed $T>0$, almost surely one has:
$$\exists N_0\geqslant 0:\forall N\geqslant N_0,\quad \forall u\in [0,T],\quad \frac{S^N_i(u)^-}{X^N_i(u)}\geqslant A_{Q_i\left(\int_{0}^u \beta I'_i(s)\mathrm{d}s\right)}$$
since $\int_{0}^u \beta I'_i(s)\mathrm{d}s\leqslant \int_{0}^T \beta I'_i(s)\mathrm{d}s$ is almost surely finite and does not depend on $N$, which yields the first part of the Theorem.\pe

To see that (\ref{conv_result}) holds, write that with probability $1$,
$$\forall T\geqslant 0,\quad \exists N_0\in \N^* ~:~ \forall N\geqslant N_0,\forall u\in [0,T],\quad I^N(u)=I'(u)$$
so 
$$\exists N_0\in \N^* ~:~ \forall N\geqslant N_0,\forall u\in [0,\tau'],\quad I^N(u)=I'(u) \text{ (so }\tau^N=\tau'\text{ and }Z^N=Z'\text{)}$$
almost surely on the part of the sample space where $\tau'<+\infty$ and $Z'<+\infty$ (note that $N_0$ is random since it intrinsically depends on the value of $\tau'$). On the other part of the sample space, $\tau'=+\infty$ so $Z'=+\infty$ (recall that Assumption \ref{assump_sub} entails that deaths occur at a positive rate as long as the total population is not zero, so infinitely many new births for the branching process $I'$ have to occur if we want $\tau'$ to be infinite) and almost surely,
$$\forall T\geqslant 0, \quad \exists N_0\in \N^* ~:~ \forall N\geqslant N_0,\quad I^N(T)=I'(T)>0 \text{\quad (so }\tau^N>T\text{) and } Z^N_T=Z'_T $$
where $Z^N_T$ and $Z'_T$ are the number of upward jumps of $\Vert I\Vert_1$ and $\Vert I'\Vert_1$ before time $T$. This entails
$$\lim_{N\to +\infty}\tau^N=+\infty \quad \text{ and } \lim_{N\to +\infty}Z^N=Z'=+\infty$$
almost surely on this part of the sample space since both $Z_T^N$ and $Z_T'$ are increasing in $T$.\pe \hfill $\square$

\indent Following traditional terminology (see Chapter 4 of \cite{AB}), we will say that the epidemic undergoes a \textit{minor outbreak} when the branching process $I'$ goes extinct and a \textit{major outbreak} if it does not. As Theorem \ref{thm_approx_sir_bp} indicates, the epidemic's extinction time and total size are of the same order of those of the branching process in the minor outbreak case and tend to infinity with $N$ in the major outbreak case. The next section is devoted to computing the probability that a major outbreak occurs.

\section{The basic reproduction number and major outbreak probability}\label{sec_major_outbreak}
\sectionmark{$R_0$ and $\P(Z'=+\infty)$}

\indent \indent A classical quantity of interest when trying to evaluate the major outbreak probability in a SIR model is the basic reproduction number $R_0$ that roughly represents the mean number of susceptibles an average infective can contaminate in the early stages of the epidemic within an otherwise initially fully susceptible population. In the one-dimensional case, standard branching process theory shows that major outbreaks may occur with positive probability if and only if $R_0>1$ (Chapter III of \cite{AthNey}, Chapter 4 of \cite{AB}), which is also the condition for the epidemic not to go extinct in the corresponding deterministic model \cite{KMK,AB,Bri10}. In the multitype case, the relevant definition for $R_0$ is the dominant eigenvalue of the mean offspring matrix for the branching process (see Chapter V of \cite{AthNey} or \cite{Ath68}). Note that all sources cited here deal with branching processes with splitting at death, but it is easy to see that their results transpose to our setting. In this section, we derive a computation method for $R_0$ and the major outbreak probability, then apply it to a real network using our BDNI data subsample.

\subsection{$R_0$ and major outbreak probability computation}\label{sec_major_outbreak_th}
\subsectionmark{Computation}

\indent \indent We start with the following Proposition on the expected offspring and (sub)criticality condition for $I'$.

\begin{prn}\label{matrix_expected_off}
\textit{For any $i,j\in \{1,\ldots,n\}$, let $W_{i,j}$ denote the number of offspring in node $j$ of an individual born in node $i$ for the branching process described by $I'$. Also set $\Sigma_i=\gamma_i+d_i+\sum_{j\neq i}\theta_{i,j}$ for any $i\in \{1,\ldots,n\}$. Let $\Theta$ be the $n\times n$ matrix defined by $\Theta_{i,j}=\theta_{i,j}$ if $i\neq j$ and $\theta_{i,i}=0$. Then the matrix $C=(\E(W_{i,j}))_{1\leqslant i,j\leqslant n}$ is given by
\begin{equation}\label{mat_exp_off}C=(\mathrm{diag}(\Sigma_1,\ldots,\Sigma_n)-\Theta)^{-1}\mathrm{diag}(\beta_1,\ldots,\beta_n).\end{equation}
Let $R_0$ denote the largest real eigenvalue of $C$. If $R_0\leqslant 1$, then the major outbreak probability is $0$. If $R_0>1$, then this probability is strictly positive and equals $1-\prod_{k=1}^n q_k^{I(0)_k}$, where $q=(q_1,\ldots,q_n)$ is the only fixed point in $[0,1)^n$ of 
$$\begin{array}{cccc}G:&[0,1]^n &\longrightarrow & [0,1]^n \\ &s=(s_1,\ldots,s_n)&\longmapsto &\left(\E\left[\prod_{j=1}^n s_j^{W_{1,j}}\right],\ldots,\E\left[\prod_{j=1}^n s_j^{W_{n,j}}\right]\right)
\end{array}.$$
Moreover, the iterated sequence $(G^{k}(s))_{k\geqslant 0}$ converges to $q$ for any $s\in [0,1)^n$.}
\end{prn}

\textsc{Proof of Proposition \ref{matrix_expected_off}} --- Although it is of paramount interest, the statement regarding the major outbreak probability and the convergence of $(G^{k}(s))_{k\geqslant 0}$ are quite classical, see Chapter 1 of \cite{Mode} or Chapter V of \cite{AthNey}.\ppe 

Our present problem is very close to that studied in \cite{Nea12}, but our strategy to compute the $\E(W_{i,j})$ has to slightly differ since we consider node-dependent death rates $d_i$ and recovery rates $\gamma_i$. We therefore introduce an additional "dead or removed" cemetery node $\partial$ and consider $I'$ as a multitype branching process in which individuals do not die anymore but move between nodes $1,\ldots,n$ at the usual rates $\theta_{i,j}$ and, while in node $i$, produce offspring at rate $\beta_i$ and jump to node $\partial$ at rate $\gamma_i+d_i$. Once in $\partial$, they cannot move from this node ($\theta_{\partial i}=0$) or have offspring ($\beta_{\partial}=0$). It is not difficult to see that our connectivity assumption entails that absorption by $\partial$ is $\P$-almost certain since at least one $d_i$ is positive. Let us now consider an individual born in node $i$. Denote by $\xi_0,\ldots,\xi_{\tau}$ his successive positions, with $\xi_0=i$ and $\xi_{\tau}=\partial$, $\tau$ being the time of absorption of the underlying Markov chain by $\partial$, and denote by $T_0,\ldots,T_{\tau}$ the times spent by the individual in these positions during the corresponding stays. The Poisson processes driving the epidemic dynamics are independent from those driving movements between nodes, so for any $j\in \ie{1}{n}$ the offspring $W_{i,j}^k$ of the individual in node $j$ during its $k$-th stay has law Poisson with parameter $\beta_jT_k$ conditionally on $(\xi_k=j)$. We may thus write:
\begin{align*}\E\left(W_{i,j}\right)&=\E\left(\sum_{k=0}^{+\infty}W_{i,j}^k\mathrm{1}_{\xi_k=j}\right)=\sum_{k=0}^{+\infty}\E\left(\E\left(W_{i,j}^k\mid \xi_k=j\right)\mathrm{1}_{\xi_k=j}\right)&\\
&=\sum_{k=0}^{+\infty}\beta_j\E\left(T_k\mathrm{1}_{\xi_k=j}\mid \xi_0=i\right)=\beta_j\sum_{k=0}^{+\infty}\E\left(\mathrm{1}_{\xi_k=j}\E\left(T_k\mid \xi_k=j, \xi_0=i\right)\right).&
\end{align*}
Now the $T_k$ follow an exponential distribution with mean $\Sigma_j^{-1}$ conditionally on $(\xi_k=j, \xi_0=i)$ so
$$\E\left(W_{i,j}\right)=\frac{\beta_j}{\Sigma_j}\E\left(\sum_{k=0}^{+\infty}\mathrm{1}_{\xi_k=j}\mid \xi_0=i\right)$$
and $\E\left(\sum_{k=0}^{+\infty}\mathrm{1}_{\xi_k=j}\mid \xi_0=i\right)=\E\left(\sum_{k=0}^{\tau-1}\mathrm{1}_{\xi_k=j}\mid \xi_0=i\right)$ is the expected number of visits of $j$ by the underlying Markov chain before absorption starting from $i$, which is known (see for instance \cite{KS60}, Chapter 3) to be $\left(\mathrm{I}_n-T\right)^{-1}_{i,j}$ where $T$ is the matrix defined by $$T_{i,j}=\begin{cases}0 \text{ if }i=j\\ \frac{\theta_{i,j}}{\Sigma_i}\text{ if }i\neq j\end{cases}$$
so
$$\E(W_{i,j})=\frac{\beta_j}{\Sigma_j}(I_n-T)^{-1}_{i,j}=\beta_j(\text{diag}(\Sigma_1,\ldots,\Sigma_n)-\Theta)^{-1}_{i,j}$$
and (\ref{mat_exp_off}) follows.
\hfill $\square$\pe

\textit{Remark :} Computing the infinitesimal generator of the mean matrix semigroup of $I'$ and using Theorem 1 from \cite{Ath68} (that, again, deals with a branching process with splitting at death but may be transposed easily to our setting) yields another necessary and sufficient condition for $I'$ to go extinct with probability $1$. The probability that a major outbreak occurs is positive if and only if the maximal real eigenvalue $\lambda_1$ of the infinitesimal generator of the mean matrix semigroup of $I'$ defined by
$$M_{i,j}=\begin{cases}\beta_i-\Sigma_i\text{ if }i=j\\ \theta_{i,j}\text{ if }i\neq j\end{cases}$$
is positive. Contrary to $R_0$, $\lambda_1$ cannot be interpreted in a straightforward biological way. Yet, it yields valuable information on the behavior of $I'$. In particular, if $\lambda_1>0$, $\lambda_1$ approximates for large $N$ the early exponential growth rate of $(I_t)_{t\geqslant 0}$ in the case where the epidemic takes off (see \cite{WL}).\pe


\indent \indent Proposition \ref{matrix_expected_off} suggests that a numerical computation of the major outbreak probability can be performed iteratively if $G$ is known. \cite{Nea12} puts forward a procedure to compute $G$ numerically when the $\gamma_i+d_i$ are identical across nodes. However, in our setting the duration of the infectious period is not independent from the infective's trajectory, hampering the integration of the conditional expectancy obtained by Theorem 1 of \cite{DM68} (see equations (11) and (12) of \cite{Nea12}). The alternative method we will now present strongly relies on the Markov property. It may be generalized to Gamma-distributed infectious periods by introducing a stage-based structure as in \cite{AB00}. Yet, it does not require diagonalisability assumptions and involves less complex computations than those of \cite{Nea12} in the cases where the latter applies.

\begin{prn}[Computation of the MGF of the $W_{i,j}$]\label{thm_calcul_proba}
\textit{Let $\Theta$ be the matrix defined in Proposition \ref{matrix_expected_off} and $\omega=d+\gamma$ the $n$-dimensional vector defined by
$$\omega_i=d_i+\gamma_i.$$
For any $s=(s_1,\ldots,s_n)$ and any $i\in \{1,\ldots,n\}$, \begin{equation}\label{eq_fonc_gen_W}G(s)=\left(\mathrm{diag}\left(\lambda_1(s),\ldots,\lambda_n(s)\right)-\Theta\right)^{-1}\omega\end{equation}
where $\lambda_i(s)=(1-s_i)\beta_i+\Sigma_i$.}
\end{prn}

\textsc{Proof of Proposition \ref{thm_calcul_proba}} --- In order for our proof to be relatively self-contained, let us first reproduce the derivation of Equation (9) from \cite{Nea12}. If $Q_{i,j}$ denotes the random time spent in node $j$ by an individual born in node $i$ (for the infectious branching process, that is, for $I'$), then the $W_{i,j}$ are independent conditionally on the $Q_{i,j}$ and $W_{i,j}$ has law Poisson with mean $\beta_j Q_{i,j}$ conditionally on $Q_{i,j}$, so for any $s\in [0,1]^n$:
\begin{align}\label{eq_9_neal}G_i(s)&:=\E\left(\prod_{j=1}^ns_j^{W_{i,j}}\right)\notag&\\
&=\E\left[\E\left(\prod_{j=1}^ns_j^{W_{i,j}}\mid Q_{i,1},\ldots,Q_{i,n}\right)\right]\notag&\\
&=\E\left[\prod_{j=1}^n\E\left(s_j^{W_{i,j}}\mid Q_{i,j}\right)\right]\notag&\\
&=\E\left[\prod_{j=1}^n \exp\left(-(1-s_j)\beta_j Q_{i,j}\right)\right],&
\end{align}
which is \cite{Nea12}'s Equation (9). Now let $(a_1,\ldots,a_n)\in \R_+^n$ and consider $I'$ as a $(n+1)$-type branching process with immortal individuals as in the proof of Proposition \ref{matrix_expected_off}. Using the law of total probability and the strong Markov property at the first jump time of an individual born in node $i$, we get:
\begin{flalign}\label{eq_fonc_gen_Q}
&\E\left(e^{-\sum_{j=1}^n a_j Q_{i,j}}\right)&&\notag\\&=\E\left(e^{-a_1T_0}\mid \xi_0=i\right)\left(\sum_{k=1}^n \P(\xi_1=k\mid \xi_0=i)\E\left(e^{-\sum_{j=1}^n a_j Q_{k,j}}\right)+\P(\xi_1=\partial\mid \xi_0=i)\right)\notag&&\\
&=\frac{\Sigma_i}{a_i+\Sigma_i}\left(\sum_{k=1}^n\frac{\theta_{i,k}}{\Sigma_i}\E\left(e^{-\sum_{j=1}^n a_j Q_{k,j}}\right)+\frac{d_i+\gamma_i}{\Sigma_i}\right)&\notag\\
&=\sum_{k=1}^n\frac{\theta_{i,k}}{a_i+\Sigma_i}\E\left(e^{-\sum_{j=1}^n a_j Q_{k,j}}\right)+\frac{d_i+\gamma_i}{a_i+\Sigma_i},&
\end{flalign}
since $T_1$ is exponentially distributed with mean $\Sigma_i^{-1}$ for an individual born in node $i$. Combining (\ref{eq_9_neal}) and (\ref{eq_fonc_gen_Q}) shows that $G(s)$ is such that
$$\left(\mathrm{diag}\left(\lambda_1(s),\ldots,\lambda_n(s)\right)-\Theta\right)G(s)=\omega,$$
but $\mathrm{diag}\left(\lambda_1(s),\ldots,\lambda_n(s)\right)-\Theta$ is a diagonally dominant matrix that can be proved to be invertible using the connectivity assumption, which ends the proof of Proposition \ref{thm_calcul_proba}.\hfill $\square$\pe

Remark: A similar conditioning argument can be used to derive (\ref{mat_exp_off}) in the proof of Proposition 4 by writing that
$$\E(Q_{i,j})=\sum_{k=1}^n \frac{\theta_{i,k}}{\Sigma_i}\E(Q_{k,j})+\frac{1}{\Sigma_i}$$
and using the fact that $\E(W_{i,j})=\beta_j \E(Q_{i,j})$ for all $i,j$.

\subsection{Numerical application}\label{sec_major_outbreak_num}

\indent \indent We now illustrate the previous theorem by computing the major outbreak probability in a metapopulation of cattle where holdings are linked by trade movements. We use a toy example where epidemiological parameters are set to values close to those of Foot-and-Mouth disease (FMD), and demographic parameters are calibrated using our cattle trade data subsample. FMD is a viral multi-species disease affecting livestock, highly infectious and easy to spread through close contacts and aerosol propagation, which leads to fast dynamics. Cases of FMD usually entail trade bans, the immediate culling of all animals in detected infected herds and ring culling (see \cite{Kee} for a review of FMD models). Yet, we use FMD here as a toy example in order to compute $R_0$ and $p$. We only account for animal trade induced transmission of the disease between herds, thus neglecting other pathways, and we do not consider any control measure.

Demographic parameters $b_i$ and $d_i$ are set to values computed using the ratios of births or deaths and the average node population (i.e. the number of animals in each holding) over the year, while the $B_i$ are proxied by the total amount of entries from outside the metapopulation considered (see data description at the beginning of Section \ref{sec_data_model});. Note that the $b_i$ equal zero for $i$ corresponding to markets or assembly centers since no birth occurs in such nodes. The same calibration method holds for setting the $\theta_{i,j}$ coefficients (using the observed flows between each pair of nodes $(i,j)$), except that we set parameters with null estimates to an arbitrarily small value (namely $10^{-6}$ years$^{-1}$) in order to make sure that $A$ is irreducible. We finally use the estimates for $\beta_i$ and $\gamma_i$ given in \cite{BdR}, that is, $\beta_i\approx .67$ days$^{-1}$ and $\gamma_i\approx\frac{1}{5.5}$ days$^{-1}$ for any $i$.\pe

\begin{figure}[h!]
\begin{center}\includegraphics[scale=.55]{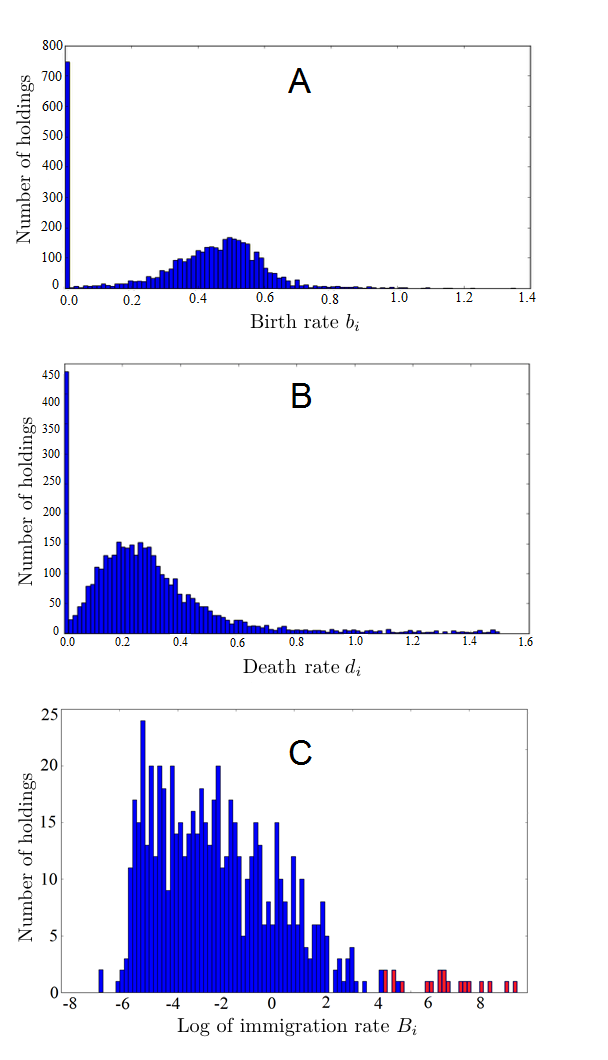}
\caption{\label{fig_histo_bidi}Values of $b_i$ (panel A), $d_i$ (panel B) and $B_i$ (panel C) in animals per year, computed using the 2015 Finistère database on cattle movements. For representability purposes, histograms on panels A and B have been truncated at $b_i=1.4$ (so 4179 holdings out of 4183 are represented), $d_i=1.5$ (3989 holdings). Panel C only accounts for the 633 nodes with non-null $B_i$. Nodes with $b_i=B_i=0$ (resp. $d_i=0$) did not exhibit any population inflow (resp. outflow) over the year. Assembly centers and markets exhibit null $b_i$, high $d_i$ not appearing on Panel B (all above 100), and high $B_i$ corresponding to observations in red in Panel C.}
\end{center}
\end{figure}

The histograms of computed values for $b_i$, $d_i$ and $B_i$ are displayed in Figure \ref{fig_histo_bidi}. Computed $b_i$ have mean .388 and standard deviation .423. Computed $d_i$ have mean 15.347 and standard deviation 407.768, but these values fall at .643 and 5.703 respectively when excluding the 20 operator nodes (that is, assembly centers and markets). Finally, computed $B_i$ have mean 1.824 and standard deviation 44.401 (respectively .117 and 2.584 without operators).\pe

\begin{figure}[!h]
\begin{center}
\includegraphics[scale=.45]{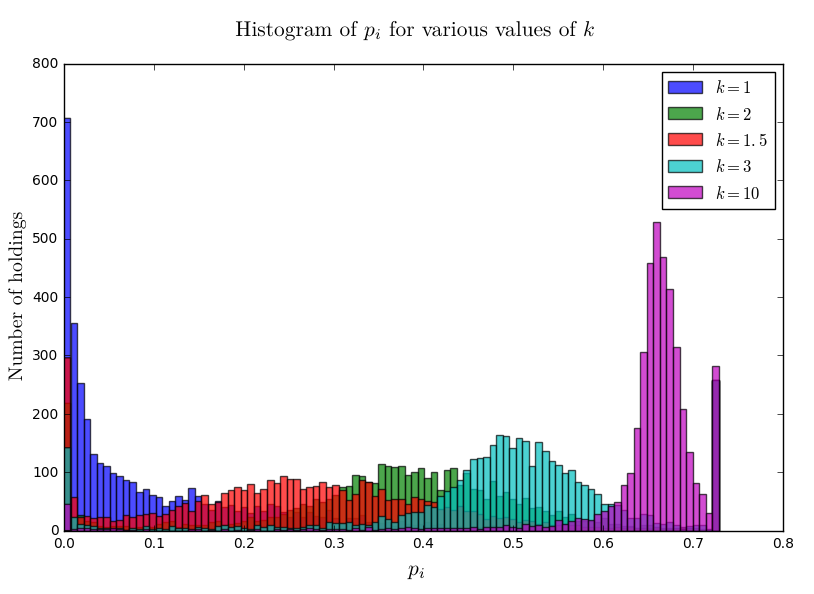}
\caption{\label{fig_probamajoroutbreak}Histograms of $p_i$ probabilities for various values of the epidemic acceleration parameter $k$. All demographic parameters are set using the 2015 Finistère database and the $\beta_i=\beta$ and $\gamma_i=\gamma$ correspond to the estimations of \cite{BdR} for foot-and-mouth disease.}
\end{center}
\end{figure}

\begin{figure}[!h]
\begin{center}
\includegraphics[scale=.33]{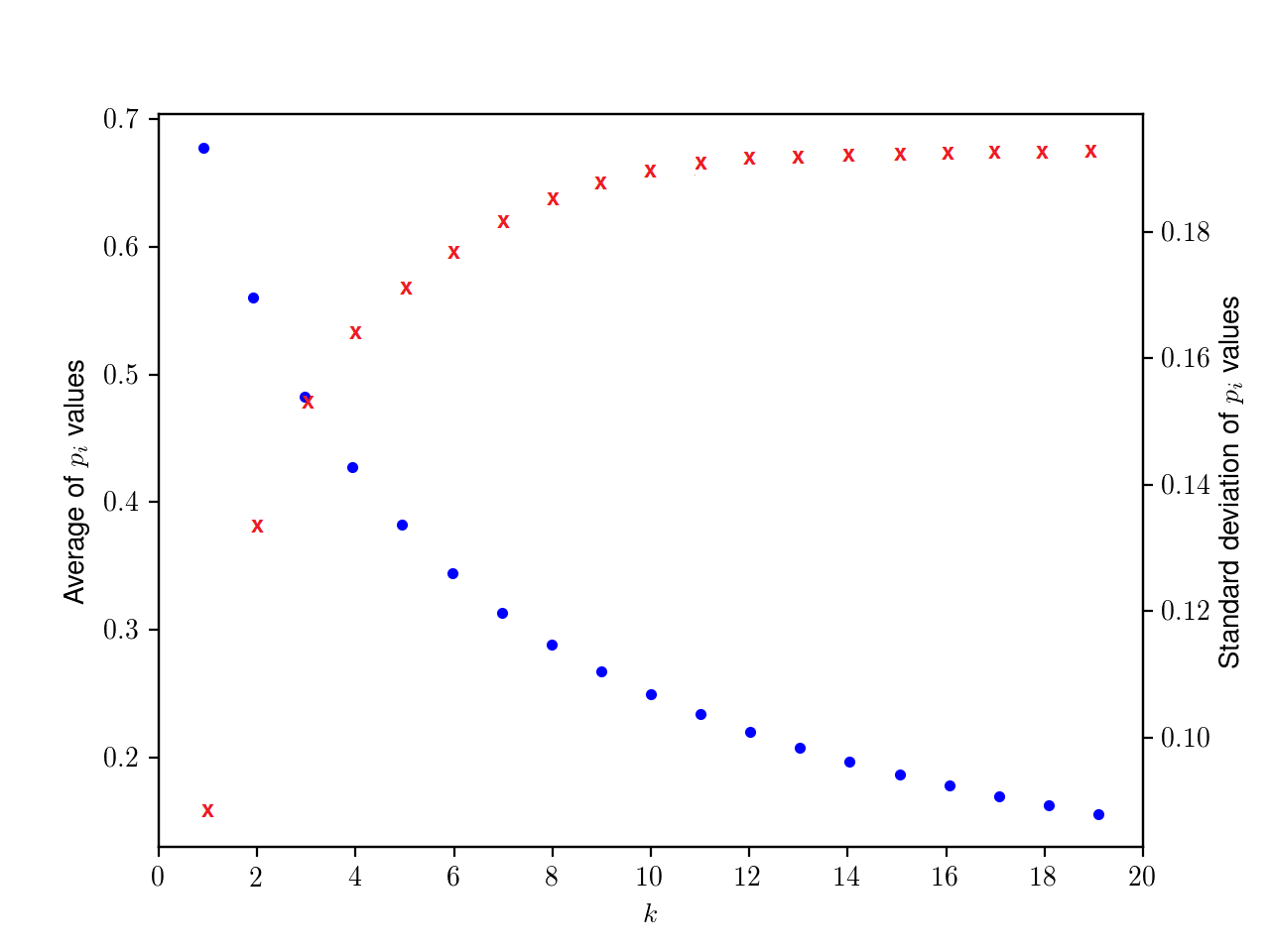}
\caption{\label{fig_meanstdpi}Mean (blue points) and standard deviation (red crosses) of $p_i$ probabilities for various values of the epidemic acceleration parameter $k$. All demographic parameters are set using the 2015 Finistère database and the $\beta_i=\beta$ and $\gamma_i=\gamma$ correspond to the estimations of \cite{BdR} for foot-and-mouth disease.}
\end{center}
\end{figure}

Iterating $G$ yields an approximation for the major outbreak probability $p:=1-q$ having mean .150, standard deviation .196 and ranging from $0$ to $.729$ over the set of nodes. Virtually accelerating the course of the epidemic by multiplying both $\beta$ and $\gamma$ by some factor $k$ (which increases the transmission rate and decreases the infection period proportionally with $k$) yields a higher mean value and a lower variance for the $p_i$, as we illustrate in Figures \ref{fig_probamajoroutbreak} and \ref{fig_meanstdpi}. Additional investigation shows that lower $p_i$ are associated with epidemics starting in nodes with very high removal rates $d_i+\gamma$ such as assembly centers, or strong transfer rates to nodes with high death rates. More generally, they are linked to nodes that are at the origin of heavily weighted paths to exit the system. As $k$ grows, the role of inter-nodal transfers decreases. Discrepancies between the $p_i$ for high $k$ are mostly due to the diversity of death rates in various nodes. Most $p_i$ related to farms get closer to the maximal value .729 corresponding to the major outbreak probability within an isolated node with $d_i=0$. This accounts for the fact that for high values of $k$, outwards movements from such nodes are on a slower time scale than the inner epidemic dynamics. For all values of $k$, the $R_0$ estimate is $3.6849$, very close to the basic reproduction number $\frac{\beta}{\gamma}\approx 3.6850$ of a closed, homogeneously mixing SIR model, which can be explained by the small magnitude of $\Theta$ as compared to $\text{diag}(\Sigma_1,\ldots,\Sigma_2)$ (see the expression for $C$ in Proposition \ref{matrix_expected_off}).

\section{The extinction time and total size of the epidemic}
\label{sec_total_size}

\indent \indent Other quantities of interest in the study of an epidemic process are its extinction time and its total size, that is, the total number of individuals infected during the course of the epidemic. In this section, we exhibit a lower bound for the maximal fraction of the population infected at a given time and use it to derive an exponential lower bound for the epidemic extinction time and total size in the case where there exists a stable endemic equilibrium for some associated dynamical system.\pe 

\indent \indent Theorem \ref{thm_approx_sir_bp} states that the total size $Z^N$ of the epidemic almost surely goes to infinity with $N$ on the part of the sample space where the approximating branching process $I'$ explodes. A standard result of the unidimensional SIR model without demography states that $Z^N$ satisfies a central limit theorem conditionally on a major outbreak occurring. In other words, major outbreaks are characterized by a positive fraction of the population being affected by the epidemic at some point in time (\cite{Sca90,AB}, and see \cite{BC93} for a multidimensional generalization) when individuals cannot enter or leave the system. The following proposition states that in the case of a major outbreak, the maximal number of infectives during the course of the epidemic is \textsl{at least} equal to a fraction of $N$ with high probability as $N$ goes to infinity when the population process starts at its equilibrium value $z^*:=-A^{-1}B$ (see Section \ref{subsec_pop_mod}). Its proof is postponed to the Appendix.

\begin{prn}\label{lower_bound_final_size}
\textit{Assume that $\tilde{x}_0=z^*$. Then
$$\lim_{\varepsilon\to 0}\sup_{N\geqslant 1}\P\left(\max_{t\in \R_+}\Vert I^N(t)\Vert_1\leqslant \varepsilon \Vert z^*\Vert_1 N, Z'=+\infty\right)=0.$$}
\end{prn}

\pe A trivial corollary to this result is that a similar lower bound also holds for the total size of the epidemic. Although this type of bound is the general rule in closed population models, it happens to be quite bad in cases where susceptible population renewal through demographic mechanisms is strong enough to entail endemicity (that is, long-term persistence of the epidemic over a given threshold). This latter phenomenon has been thoroughly investigated for deterministic systems in the single-population case (see the deterministic approximations in \cite{AB00,Nas,VH95,VH97}); in our stochastic framework, we may reasonably expect an increase of $N$ to affect both the typical infective population during the course of the epidemic and the time scale of the endemic period, thus resulting in a more than proportional total size response.\pe

Let us consider the dynamical system on $E=\{(s,i,r)\in \R_+^{3n}, i\neq 0\}$ defined by:
\begin{equation}\label{dynsys}\tag{$\mathcal{S}$}\left.\begin{array}{ccl}
\dot{s}_k&=&B_k+b_k(s_k+i_k+r_k)-d_k s_k +\sum_{j\neq k}\theta_{j,k}s_j-\sum_{j\neq k}\theta_{k,j}s_k- \beta_k i_k \frac{s_k}{z^*_k}\ppe\\
\dot{i}_k&=&\beta_k i_k\frac{s_k}{z^*_k} - d_k i_k - \gamma_k i_k+\sum_{j\neq k}\theta_{j,k}i_j-\sum_{j\neq k}\theta_{k,j}i_k\ppe\\
\dot{r}_k&=&\gamma_k i_k - d_k r_k+\sum_{j\neq k}\theta_{j,k}r_j-\sum_{j\neq k}\theta_{k,j}r_k
\end{array}\right\}\end{equation}
for all $k\in \ie{1}{n}$.\pe

The following result shows that the existence of a globally attractive endemic equilibrium for ($\mathcal{S}$) yields an exponential lower bound for the extinction time and the total size of the epidemic in the major outbreak scenario.\pe 

\begin{thth}\label{thm_syst_dyn_exp}
\textit{Assume that $x_0=z^*$ and that ($\mathcal{S}$) admits a globally asymptotically stable endemic equilibrium $(s^*,i^*,r^*)\in E$.
Then:
\begin{equation}\label{endem_exp_tps}
\exists a_1>0 : \forall \alpha>0, ~\lim_{N\to +\infty}~\P\left(\tau^N \leqslant e^{(a_1-\alpha) N}, Z'=+\infty\right)=0.\end{equation}
and
\begin{equation}\label{endem_exp}
\exists a_1>0 : \forall \alpha>0, ~\lim_{N\to +\infty}~\P\left(Z^N \leqslant e^{(a_1-\alpha) N}, Z'=+\infty\right)=0.\end{equation}}
\end{thth}

\textit{Remark :} if the assumption of Theorem \ref{thm_syst_dyn_exp} is met, a quick look at ($\mathcal{S}$) shows that $i^*$ has only strictly positive coordinates. Then, as we will see in the following proof, $a_1$ is the exit cost from $\mathcal{B}_{\infty}\left((s^*,i^*,r^*),\min_j i^*_j\right)$ for ($\mathcal{S}$), given (see Equation (2.3) of Chapter 5 of \cite{FW}, \cite{ParBrice} or \cite{KP}) by

\begin{equation*}a_1=\inf_{y\in \partial \mathcal{B}_{\infty}((s^*,i^*,r^*),\min_j i^*_j)}~\inf_{\phi}~\int_{T_1}^{T_2}L(\phi_t,\dot{\phi}_t)\mathrm{d}t,\end{equation*}
where the second infimum is taken over the set of absolutely continuous functions $\phi$ on some $[T_1,T_2]$ (with $-\infty\leqslant T_1<T_2\leqslant +\infty$) such that $\phi(T_1)=(s^*,i^*,r^*)$ and $\phi(T_2)=y$, and $L$ is defined by 
\begin{align*}&L((s,i,r),\beta):=\max_{u\in \R^{3n}}\left[\beta\cdot u - \sum_{j=1}^n (e^{u_j}-1)(B_j+b_j(s_j+i_j+r_j))\right.&&\\
&\quad-\sum_{j=1}^n \left[(e^{-u_j}-1)d_j s_j+(e^{-u_{n+j}}-1)d_j i_j+(e^{-u_{2n+j}}-1)d_j r_j\right]&&\\
&\quad- \sum_{j\neq k} \left[(e^{u_k-u_j}-1)\theta_{j,k}s_j+(e^{u_{n+k}-u_{n+j}}-1)\theta_{j,k}i_j+(e^{u_{2n+k}-u_{2n+j}}-1)\theta_{j,k}r_j\right]&&\\
&\quad-\sum_{j=1}^n \left[(e^{-u_j}-1)d_j s_j+(e^{-u_{n+j}}-1)d_j i_j+(e^{-u_{2n+j}}-1)d_j r_j\right]&&\\
&\quad-\sum_{j=1}^n (e^{u_{n+j}-u_j}-1)\beta_j\frac{s_ji_j}{s_j+i_j+r_j}\left.-\sum_{j=1}^n (e^{u_{2n+j}-u_{n+j}}-1)\gamma_ji_j\right]&&\end{align*}
for all $(s,i,r)\in \R_+^{3n}$ and $\beta\in \R^{3n}$.\pe

When $n\geqslant 2$, the proof of Theorem 4.1 from \cite{LS} adapts and shows that if $R_0>1$ and if there exists $\lambda>0$ such that $s^*=\lambda i^*$, then the endemic equilibrium is globally asymptotically stable and the conclusion of Theorem \ref{thm_syst_dyn_exp} is true.\pe

\textsc{Proof of Theorem \ref{thm_syst_dyn_exp}} --- The conclusion of Theorem \ref{thm_syst_dyn_exp} is obviously true if $R_0\leqslant 1$, so we now assume that $R_0>1$. Let $\delta=\min_j i^*_j>0$ and let $\delta'>0$. Let $\varepsilon>0$ be such that setting
$$T_{\varepsilon}=\inf~\{t\geqslant 0 : \Vert I^N(t)\Vert_1\geqslant \varepsilon N\}$$
then
\begin{equation}\label{prem_liminf}\liminf_{N\to +\infty}~\P\left(T_{\varepsilon}<+\infty\right)>\P(Z'=+\infty)-\delta'.\end{equation}
Such a $\varepsilon$ exists according to Proposition \ref{lower_bound_final_size}. Next define the flow $\Phi:E\times \R_+\to E$ associated with the dynamical system. There exists $M>1$ such that $\Vert z_t\Vert_1$ only takes values below $M\Vert z^*\Vert_1$ for $t\geqslant 0$ whenever $z_0$ is close enough from $z^*$. There also exists $T>0$ such that: 
\begin{multline}\label{temps_entre_boule}\forall (s,i,r)\in E : \Vert i\Vert_1\geqslant\varepsilon \text{ and } \Vert s+i+r \Vert_1\leqslant M\Vert z^* \Vert_1,\\ \Phi((s,i,r),T)\in \mathcal{B}_{\infty}\left((s^*,i^*,r^*),\frac{\delta}{2}\right),\qquad \qquad \qquad\end{multline}
since $(s,i,r)\mapsto \inf~\left\{t\geqslant 0 : \Phi((s,i,r),t)\in \mathcal{B}_{\infty}\left((s^*,i^*,r^*),\frac{\delta}{2}\right)\right\}$ is upper semi-continuous (and then upper bounded) on $\{(s,i,r) \in E : \Vert i\Vert_1\geqslant\varepsilon, \Vert s+i+r \Vert_1\leqslant M\Vert z^* \Vert_1\}$ that is a compact set. As a result:
\begin{equation}\label{entre_boule}\Phi\left(\frac{(S^N,I^N,R^N)}{N}({T_{\varepsilon}}),T\right)\in \mathcal{B}_{\infty}\left((s^*,i^*,r^*),\frac{\delta}{2}\right)\end{equation}
almost surely conditionally on $\left(T_{\varepsilon}<+\infty, \Vert X^N({T_{\varepsilon}})\Vert_1\leqslant M\Vert z^*
\Vert_1 N\right)$. Mimicking the proof of Theorem 2.1 from Chapter 11 of \cite{EK}, we get that
$$\sup_{t\in [0,T]}\left\Vert \frac{(S^N,I^N,R^N)}{N}({T_{\varepsilon}+t})-\Phi\left(\frac{(S^N,I^N,R^N)}{N}({T_{\varepsilon}}),t\right)\right\Vert_{\infty}< \frac{\delta}{2}$$ 
and
$$\Vert X^N({T_{\varepsilon}})\Vert_1\leqslant M\Vert z^*
\Vert_1 N$$
with probability going to $1$ conditionally on $(T_{\varepsilon}<+\infty)$ when $N\to +\infty$.
This, (\ref{prem_liminf}), (\ref{entre_boule}) and the definition of $T$ show that 
\begin{align}\label{lower_bd_hit}&\liminf_{N\to +\infty}~ \P\left(T_{\varepsilon}<+\infty,\exists t\geqslant 0 : \frac{(S^N(t),I^N(t),R^N(t))}{N}\in \mathcal{B}_{\infty}\left((s^*,i^*,r^*),\delta\right)\right)&\notag
\\&\quad>\P(Z'=+\infty)-\delta',&\end{align}
so $\frac{(S^N,I^N,R^N)}{N}$ hits $\mathcal{B}_{\infty}\left((s^*,i^*,r^*),\delta\right)$ with probability at least $\P(Z'=+\infty)-\delta'$ for $N$ large enough. Using Theorem 6 of \cite{ParBrice} just as in the proof of Proposition \ref{thm_fw_pop} along with the Markov property yields $a_1>0$ such that for all $\alpha>0$,
\begin{align}\label{temps_passe_sur_seuil}&\liminf_{N\to +\infty}~ \P\left(\exists T\geqslant 0, \forall t\in [T,T+e^{(a_1-\alpha)N}], \forall j\in \{1,\ldots,n\}, I^N_j(t)>(i^*_j-\delta) N\right)\notag
\\&\quad>\P(Z'=+\infty)-\delta',&\end{align}
and (\ref{endem_exp_tps}) follows.
\pe 
Now set $\alpha>0$ and take $\alpha'\in (0,\alpha)$. Then, conditionally on the event $$\left(\exists T\geqslant 0, \forall t\in [T,T+e^{(a_1-\alpha')N}], \forall j\in \{1,\ldots,n\}, I^N_j(t)>(i^*_j-\delta) N\right),$$ the total number of infected individuals recovering or dying during the course of the epidemic stochastically dominates the value at time $e^{(a_1-\alpha')N}$ of a homogeneous Poisson counting process $(Q^N(t))_{t\geqslant 0}$ with intensity $\lambda_N:=\min_{j}(\gamma_j+d_j)(i^*_j-\delta) N$. Yet the former number is also lower than $Z^N$ with probability $1$ since all infected individuals eventually have to die or recover, so we finally get:
\begin{align*}\liminf_{N\to +\infty}\, \P\left(Z^N\geqslant e^{(a_1-\alpha)N}\right)&\geqslant \liminf_{N\to +\infty}\, \P\left(Q^N_{e^{(a_1-\alpha)N}}\geqslant e^{(a_1-\alpha')N}\right)\left(\P(Z'=+\infty)-\delta'\right)&
\end{align*}
using (\ref{temps_passe_sur_seuil}), which yields
$$\liminf_{N\to +\infty}~\P\left(Z^N\geqslant e^{(a_1-\alpha)N}\right) \geqslant \P(Z'=+\infty)-\delta',$$
because the law of $Q^N\left(e^{(a_1-\alpha')N}\right)$ is Poisson with mean $\lambda_N e^{(a_1-\alpha')N}$, hence (\ref{endem_exp}). 
\hfill $\square$

\pe

Considering the proof of Proposition \ref{lower_bound_final_size}, we can see that the time needed for $I^N$ to go above a given fraction of $N$ is of order $\log(N)$ on the event $Z'=+\infty$. The proof of Theorem \ref{thm_syst_dyn_exp} shows that the subsequent convergence time of $\frac{(S^N,I^N,R^N)}{N}$ towards a given ball centered on the endemic equilibrium is upper bounded by some deterministic constant $T$ with high probability, while the results from \cite{FW} state that the time needed for the scaled process to leave the ball is of order $e^{a_1 N}$. This shows that for large $N$, the epidemic undergoing a major outbreak spends most of its time in its endemic phase where the scaled process lies close to the endemic equilibrium. As a result, the lower bound for $Z^N$ we found in Theorem \ref{thm_syst_dyn_exp} appears to be of the right order.

\pe

When $n=1$, direct calculations yield a simple necessary and sufficient condition for the existence of a globally stable endemic equilibrium. This leads to the following corollary, the proof of which is given in the Appendix.

\begin{coe}\label{thm_syst_dyn_exp_dim1}
\textit{Assume that $n=1$ and $\beta>d+\gamma$, so the major outbreak probability is positive (we omit the $\cdot_1$ subscripts in the parameters). Then any solution $(s,i,r)$ of
\begin{equation}\label{simili_lv}\begin{cases}
\dot{s}=B+b(s+i+r)-d s-\beta i\frac{s}{s+i+r} \\
\dot{i}=\beta i \frac{s}{s+i+r} - di-\gamma i\\
\dot{r}=\gamma i - dr
\end{cases}\end{equation}
converges to the endemic steady state
$$(s^*,i^*,r^*)=\begin{pmatrix}\frac{d+\gamma}{\beta}\frac{B}{d-b},d\frac{B}{d-b}\left(\frac{1}{d+\gamma}-\frac{1}{\beta}\right),\gamma\frac{B}{d-b}\left(\frac{1}{d+\gamma}-\frac{1}{\beta}\right)\end{pmatrix}$$
and the conclusions of Theorem \ref{thm_syst_dyn_exp} hold (recall that Assumption \ref{assump_sub} implies $d>b$).}
\end{coe}

\section{Conclusion}

\indent \indent We defined a multitype, stochastic SIR dynamical epidemic model on a strongly connected graph. Using a branching approximation, we defined minor and major epidemic outbreaks and gave a necessary and sufficient condition for major outbreaks to occur, along with a computational method for the probability of such events when the condition if fulfilled. Our main result consists of an exponential lower bound for the extinction time and the total size of the epidemic in the stable endemic case when a major outbreak occurs, improving on the usual results for demography-free dynamics. 

Although Theorem \ref{thm_syst_dyn_exp} gives a rather good lower bound for the size of the epidemic in the major outbreak case, we do not know much about its distribution yet and we are still investigating on the transposition of results from \cite{Sca90,BC93} to our open, multinodal setting.

Another issue we plan to discuss on in the foreseeable future is the existence of a quasi-stationary distribution for the epidemic process, that is, of an asymptotic distribution conditionally on non-extinction \cite{VJ,VDP13}. Proving that such a distribution exists seems challenging even in the open uninodal cases considered by \cite{AB00} and \cite{Nas}, and usual criteria do not apply.

Both questions might benefit from considering a diffusive scaling limit of our model \cite{VH95,VH95,Nas,LJ}, making it possible to use Fokker-Plank equations for computing the fade out probability of an epidemic after its first major outbreak, and to derive quasi-stationary approximations for the limiting diffusion.

Finally, designing an efficient simulation and estimation procedure and estimation procedure in order to calibrate the model on data will make it possible to illustrate the results of Section \ref{sec_total_size} as we did for the results of Section \ref{sec_major_outbreak_th} while avoiding prohibitive computation times.

\section*{Acknowledgements}

\indent \indent This work is part of a PhD Thesis supervised by Vincent Bansaye (CMAP, École Polytechnique) and Elisabeta Vergu (MaIAGE, INRA), whom I warmly thank for their guidance and support. It was supported by the French Research Agency within projects ANR-16-CE32-0007-01 (CADENCE) and ANR-16-CE40-0001 (ABIM), and by Chaire Modélisation Mathématique et Biodiversité Veolia-X-MNHM-FX.

\section{Appendix}

\subsection{Proof of Proposition \ref{th_stabpop}}

This statement is reminiscent of a classical result on multitype branching processes (see Chapter IV.7 of \cite{AthNey} or Theorem 4.2.2 of \cite{Jag}). However, the definition of such processes slightly differs from the one we chose here: in the classical setting, individuals do not move between nodes and only split at death between other individuals of various types. Chapter 4 of \cite{Jag} considers a unidimensional \textit{general branching process} that allows individuals to give birth at random times of their lives, and his proof could be adapted to fit our framework. One could also consider a time-sampled version of $X^N$ to retrieve a multitype Bienaymé-Galton-Watson process with immigration (see Chapter III.6 of \cite{AthNey}) or compare $X^N$ to multitype branching processes with splitting at death. Yet, we need a finer description of the return time to compact subsets of $\Z_+^n$ in order to establish not only positive recurrence but also uniform ergodicity.\pe

Proving Proposition \ref{th_stabpop} is easy when $B=0$, using that $\frac{\mathrm{d}}{\mathrm{d}t}\E(X^N(t))=A\E(X^N(t))$ and a generalized eigendecomposition of $A$ (in this particular case $\pi=\delta_0$), so we now assume that $B\neq 0$. We first show the following lemma.

\begin{lea}[Existence of a Lyapunov function for $X^N$]\label{lem_liap}
\textit{There exists $v\in \R^n$ with positive coordinates, $c>0$ and $R\geqslant 0$ such that $$v\cdot (Ax+B)<-c(v\cdot x+1)$$ for any $x\in \R_+^n$ such that $\Vert x \Vert_1\geqslant R$.}
\end{lea}
\textsc{Proof of Lemma \ref{lem_liap}} --- Recall that the transpose ${}^tA$ of $A$ is invertible because of Assumption \ref{assump_sub} and set $u=-{}^tA^{-1}B$. Quick calculations show that $u$ is the limit value of solutions of the $n$-dimensional linear ODE \begin{equation}\label{EDO}y'={}^tAy+B\end{equation} since ${}^tA$'s eigenvalues have negative real parts. Let us consider a solution of (\ref{EDO}) such that $y_0$ has positive coordinates. Writing (\ref{EDO}) as
$$\forall i\in \{1,\ldots,n\},\quad y'_i=\left[b_i-d_i-\sum_{j\neq i}\theta_{i,j}\right]y_i+\sum_{j\neq i}\theta_{i,j}y_j+B_i$$ and using that the graph with edge set $\{(i,j)\mid \theta_{i,j}>0\}$ is connected (so that all $\theta_{i,j}$ cannot be zero), we see that no $y_i$ ever reaches $0$ in finite time. As a result, $u$ has nonnegative coordinates. Similarly, if $u_i=0$ then $u_j=0$ for any $j$ such that $\theta_{i,j}>0$, hence $u=0$ by induction because of the graph connectivity, which contradicts the fact that $B\neq 0$. All components of $u$ are therefore positive. Now ${}^tAu=-B$, $B$ has nonnegative components and ${}^tA$ is invertible, so for any $x\in \R_+^n$ one may find $v$ in a neighborhood of $u$ and $C$ in a neighborhood of $B$ such that both $v$ and $C$ have positive components and such that ${}^tAv=-C$. This rewrites ${}^t v A=-{}^tC$, so 
$$v\cdot (Ax+B)=-C\cdot x+v\cdot B.$$
Defining
$$c=\frac{\min_i C_i}{2\max_i v_i}>0$$
and 
$$R=\frac{1}{\min_i v_i}\left(1+\frac{v\cdot B}{c}\right)$$
then yields the result.\pe

\textsc{Proof of Proposition \ref{th_stabpop}} --- Lemma \ref{lem_liap} shows that $f:x\mapsto 1+v\cdot x$ satisfies Condition (CD2) from \cite{MT3} with $V=f$, $C=\{x\in \Z_+^n\mid \Vert x\Vert_1\leqslant R\}$ and $d=\max_i v_i R$ since the infinitesimal generator of $\mathcal{A}$ of $X^N$ is such that $\mathcal{A}f(x)=v\cdot(Ax+B)$ for any $x\in \Z_+^n$. Moreover, it is not difficult to see that all compact sets of $\Z_+^n$ are $\delta_0$-petite for any skeleton chain of $X^N$ by considering sequences of appropriate transferts and death events --- recalling that Assumption \ref{assump_sub} implies that at least one of the $d_i$ is positive. Theorem 4.2 of \cite{MT3} thus shows that $X^N$ is positive Harris recurrent and Theorem 7.1 yields the expected result. The proposition on the first moment of $\pi$ comes from the fact that $\int (Ax+B)\mathrm{d}\pi(x)=0$ since $x\mapsto Ax+B$ is the value of the generator of $X^N$ applied to $\mathrm{Id}$.\pe \hfill$\square$

\subsection{Proof of Proposition \ref{thm_fw_pop}}

Proposition \ref{thm_fw_pop} is Theorem 6 from \cite{ParBrice} applied to a modified version of $X^N/N$ with rates vanishing outside of $\mathcal{A}=\mathcal{B}_2\left(z^*,2\varepsilon\right)\cap \R_+^n$, for instance the scaled process $\tilde{X}^N/N$ where $\tilde{X}^N$ is defined from the same Poisson processes and with the same initial condition as $X^N$ with all rates in (\ref{tableau_taux_trans}) multiplied by $$\sigma(x/N)=\mathrm{1}_{\Vert x/N-z^*\Vert_2\leqslant \varepsilon}+\mathrm{1}_{\varepsilon<\Vert x/N-z^*\Vert_2\leqslant 2\varepsilon}\left(2-\frac{\Vert x/N-z^*\Vert_2}{\varepsilon}\right).$$ The trajectories of $\tilde{X}^N/N$ are the same as those of $X^N/N$ until $\tau_{\varepsilon}^N$, so it is sufficient to apply Theorem 6 from \cite{ParBrice} to $\tilde{X}^N/N$. Note that $\mathcal{A}$ lies in the domain of attraction of $z^*$ for the dynamical system $z'=\sigma(z)Az+B$ since $A$ is negative definite. What remains to be shown is that $\alpha_0=\overline{V}:=\min_{y\in \partial \mathcal{B}_2(z^*,\varepsilon)}V(z^*,y)$ is positive, $V(z^*,\cdot)$ denoting the quasipotential of the dynamical system with respect to the Poisson perturbation (according to the terminology of \cite{FW}), defined in Section 5 of \cite{ParBrice}.\pe

It follows from $A$ having only eigenvalues with negative real parts that there exists $\varepsilon',\eta\in \left(0,\frac{\varepsilon}{2}\right]$ such that for any absolutely continuous function $\phi:\R\to \R^n$, if $\varepsilon'<\Vert \phi_t-z^*\Vert<\varepsilon$ then $\frac{d}{dt}\Vert \phi_t-z^*\Vert^2_2<0$ whenever $\Vert \dot{\phi}_t-(A\phi_t+B)\Vert_2<\eta$. Section 4 of \cite{KP} (or Equation (2.3) of Chapter 5 of \cite{FW}) now implies that
\begin{equation}\label{mino_v}\alpha_0\geqslant \inf_{y\in \partial \mathcal{B}_2(z^*,\varepsilon)}\inf_{y'\in \partial \mathcal{B}_2(z^*,\varepsilon')}\inf \int_{T_1}^{T_2}L(\phi_t,\dot{\phi}_t)\mathrm{d}t,\end{equation}
where the third infimum is taken over the set of $\{x\in \R^n : \varepsilon'\leqslant \Vert x-z^*\Vert_2 \leqslant  \varepsilon\}$-valued absolutely continuous functions $\phi$ on some $[T_1,T_2]$ (with $-\infty\leqslant T_1<T_2\leqslant +\infty$) such that $\phi(T_1)=y'$ and $\phi(T_2)=y$, and $L$ is defined by \begin{multline}\label{def_L}L(x,\beta):=\max_{u\in \R^n}\left[\beta\cdot u - \sum_i (e^{u_i}-1)(B_i+b_i x_i)\right.\\ \left.-\sum_i (e^{-u_i}-1)d_i x_i - \sum_{i\neq j} (e^{u_j-u_i}-1)\theta_{i,j}x_i\right]\end{multline}
for all $x\in \R_+^n$ and $\beta\in \R^n$. Now let $y\in \partial \mathcal{B}_2(z^*,\varepsilon)$ and $y'\in \partial \mathcal{B}_2(z^*,\varepsilon')$ and assume that $\inf \int_{T_1}^{T_2}L(\phi_t,\dot{\phi}_t)\mathrm{d}t=0$, with the infimum defined as before. For any choice of $\phi$ if follows from the definition of $\varepsilon'$ and $\eta$ that:
\begin{align*}\varepsilon^2-\varepsilon'^2&\leqslant \Vert y-z^*\Vert^2_2-\Vert y'-z^*\Vert^2_2& \\
&\leqslant \int_{t: \Vert \dot{\phi}_t-(A\phi_t+B)\Vert_2> \eta}\frac{\mathrm{d}}{\mathrm{d}t}\Vert \phi_t-z^*\Vert^2_2\mathrm{d}t&\\
&= 2\int_{t: \Vert \dot{\phi}_t-(A\phi_t+B)\Vert_2> \eta}(\phi_t-z^*)\cdot \dot{\phi}_t\mathrm{d}t&\\
&\leqslant 2\varepsilon \int_{t:\Vert \dot{\phi}_t-(A\phi_t+B)\Vert_2> \eta}\Vert \dot{\phi}_t\Vert_2\mathrm{d}t,&
 \end{align*}
by the Cauchy-Schwarz inequality, so
\begin{equation}
\int_{t:\Vert \dot{\phi}_t-(A\phi_t+B)\Vert_2>\eta}\Vert \dot{\phi}_t\Vert_2\mathrm{d}t\geqslant \frac{1}{2}\frac{\varepsilon^2-\varepsilon'^2}{\varepsilon}\geqslant \frac{3\varepsilon}{8}.\label{majo_int_surepsilon}
\end{equation}
Now there exists $\delta>0$ such that for all $x\in \mathcal{B}_2(z^*,\varepsilon)$ and all $\beta\in \R^n$,
\begin{equation}\label{mino_L_1}L(x,\beta)\geqslant \delta \left(\Vert \beta-(Ax+B)\Vert_2-\frac{\eta}{2}\right),\end{equation}
as seen by considering $u=\delta \frac{\beta-(Ax+B)}{\Vert \beta-(Ax+B) \Vert_2}$ in (\ref{def_L}) if $Ax+B\neq \beta$ and using a Taylor expansion for $\delta\approx 0$ for the function maximized in (\ref{def_L}), so
\begin{equation}\label{mino_L_2}L(x,\beta)\geqslant \delta \left(\Vert \beta\Vert_2 - \Vert Ax+B\Vert_2 - \frac{\eta}{2}\right).\end{equation}
Let $M>0$. If $\phi$ is such that 
$$\int_{T_1}^{T_2} L(\phi_t,\dot{\phi}_t)\mathrm{d}t < \frac{\delta\eta}{2M},$$
then the Lebesgue measure of $\{t\in [T_1,T_2] : \Vert \dot{\phi}_t-(A\phi_t+B)\Vert_2\geqslant \eta\}$ has to be lower than $\frac{1}{M}$ because of (\ref{mino_L_1}). For such a $\phi$:
\begin{align*}\int_{T_1}^{T_2}L(\phi_t,\dot{\phi}_t)\mathrm{d}t& \geqslant \int_{t : \Vert \dot{\phi}_t-(A\phi_t+B)\Vert_2> \eta}L(\phi_t,\dot{\phi}_t)\mathrm{d}t &\\
& \geqslant  \int_{t : \Vert \dot{\phi}_t-(A\phi_t+B)\Vert_2> \eta}\delta\left(\Vert \dot{\phi}_t\Vert_2-\Vert A\phi_t+B\Vert_2 - \frac{\eta}{2}\right)\mathrm{d}t& \\
& \geqslant \delta \left(\frac{3\varepsilon}{8}-\frac{\sup_{x\in \mathcal{B}_2(z^*,\varepsilon)}\Vert Ax+B\Vert_2+\frac{\eta}{2}}{M}\right),
\end{align*}
using (\ref{majo_int_surepsilon}), which contradicts, for $M$ large enough, the fact that $\int_{T_1}^{T_2}L(\phi_t,\dot{\phi}_t)\mathrm{d}t$ can be made arbitrarily small for some choice of $\phi$. This and (\ref{mino_v}) yield $\alpha_0>0$ since $\delta$ does not depend from the choice of $y$ and $y'$, which ends the proof.

\subsection{Proof of Proposition \ref{lower_bound_final_size}}

Proposition \ref{lower_bound_final_size} is obvious if $R_0\leqslant 1$, that is, if $\P(Z'=+\infty)=0$, so we now assume that $R_0>1$ (so $\P(Z'=+\infty)>0$). Let $\delta\in \left(0,\frac{\P(Z'=+\infty)}{2}\right)$. Let $\eta \in \left(0,\frac{R_0-1}{2R_0}\right)$ and $\alpha\in \left(\eta,\frac{R_0-1}{R_0}-\eta\right)$ be such that a branching process obtained from $I'$ by replacing the birth rates $\beta_i$ by $\beta_i\frac{1-\alpha}{1+\eta}$ survives with probability at least $\P(Z'=+\infty)-\delta$ (the existence of such a value of $\alpha$ is an easy consequence of Theorem \ref{thm_approx_sir_bp} and Proposition \ref{thm_calcul_proba}). We will show that if $\eta$ is small enough, then
\begin{equation}\label{a_montrer}\liminf_{N\to +\infty}\,\P\left(\max_{t\geqslant 0} \Vert I^N(t)\Vert_1 \geqslant (\alpha-\eta)\Vert z^*\hspace{-1pt}\Vert_{1}N\right)\geqslant \P\left(Z'=+\infty\right)-2\delta,\end{equation}
which implies Proposition \ref{lower_bound_final_size}. Set
$$\sigma^N_{\eta}=\inf\left\{t\geqslant 0 : \Vert I^N(t)+R^N(t)\Vert_1\geqslant (\alpha-\eta)\Vert z^*\hspace{-1pt}\Vert_{1} N\right\}.$$
Using the notation of Proposition \ref{thm_fw_pop} and setting $\eta'=\eta\Vert z^*\Vert_1n^{-1/2}$, until time $\sigma^N_{\eta}\wedge \tau^N_{\eta'}$ every infective in node $i$ makes infectious contacts with other individuals in node $j$ at rate at least $\beta_j\frac{1-\alpha}{1+\eta}$ since this node contains at least $(1-\alpha)\Vert z^*\hspace{-1pt}\Vert_{1}N$ susceptibles out of at most $(1+\eta)\Vert z^*\hspace{-1pt}\Vert_{1}N$ individuals. We can therefore define a $\R_+^{2n}$-valued multitype branching process $(I''(t),R''(t))_{t\geqslant 0}$ with rates given by 
$$\begin{array}{cc}
\text{Transition} & \text{Rate at state }(i,r) \\
(i,r)\to (i,r)-e^i_j & d_ji_j \\
(i,r)\to (i,r)-e^r_j & d_jr_j \\
(i,r)\to (i,r)+e^i_k-e^i_j & \theta_{j,k}i_j \\
(i,r)\to (i,r)+e^r_k-e^r_j & \theta_{j,k}r_j \\
(i,r)\to (i,r)+e^i_j & \beta_j \frac{1-\alpha}{1+\eta} \\
(i,r)\to (i,r)+e^r_j-e^i_j & \gamma_j i_j \\
\end{array}$$ 
and such that $\Vert I''(t\wedge \sigma_{\eta}^N\wedge \tau_{\eta'}^N)\Vert_1\leqslant \Vert I^N(t\wedge \sigma_{\eta}^N\wedge \tau_{\eta'}^N)\Vert_1$ and $\Vert R''(t\wedge \sigma_{\eta}^N\wedge \tau_{\eta'}^N)\Vert_1\leqslant \Vert R^N(t\wedge \sigma_{\eta}^N\wedge \tau_{\eta'}^N)\Vert_1$ for all $t\geqslant 0$ almost surely (so $(I'',R'')$ may not go to infinity outside of the event $(Z'=+\infty)$). Therefore:
\begin{align}\label{mino_par_bp}\P\left(\max_{t\geqslant 0} \Vert I(t)\Vert_1 \geqslant (\alpha-\eta)\Vert z^*\hspace{-1pt}\Vert_1N\right) & \geqslant  \P\left(\exists t<\tau''^N_{\eta}\wedge \sigma''^N_{\eta} : \Vert I''(t)\Vert_1 \geqslant (\alpha-\eta)\Vert z^*\hspace{-1pt}\Vert_1N\right)
\end{align}
with
$$\tau''^N_{\eta}=\inf\, \{t\geqslant 0 : \Vert I''(t) \Vert_1 \geqslant \eta \Vert z^*\Vert_1 N\}$$
and
$$\sigma''^N_{\eta}=\inf\, \left\{t\geqslant 0 : \Vert I''(t)+R''(t)\Vert_1\geqslant (\alpha-\eta)\Vert z^*\hspace{-1pt}\Vert_1 N\right\}.$$
Now $(I'',R'')$ is a non-explosive jump process so $\sigma''_{\eta}$ goes to infinity almost surely with $N$. Moreover, a continuous-time version of Theorem 2.1 of \cite{KS67} (derived for instance from this Theorem by sampling the $(I'',R'')$ at its jump times to obtain a discrete-time decomposable branching process) shows that there exists $\lambda \in (\R_+^*)^n$ such that:
$$\P\left(Z'=+\infty, \forall j\in \{1,\ldots,n\} : \lim_{t\to +\infty}\frac{I''_j(t)}{R''_j(t)}=\lambda_j\right)\geqslant \P\left(Z'=+\infty\right)-\delta$$
and that for small enough values of $\eta$, 
\begin{equation}\label{eps_small_enough}
\lambda_j>\frac{\eta}{\alpha-\eta}
\end{equation} for all $j\in \{1,\ldots,n\}$. Hence
$$\P\left(Z'=+\infty,\forall j\in \{1,\ldots,n\} : \lim_{N\to +\infty}\frac{I''^j(\sigma''^N_{\eta})}{R''^j(\sigma''^N_{\eta})}=\lambda_j\right)\geqslant \P\left(Z'=+\infty\right)-\delta$$
and (\ref{eps_small_enough}) and Fatou's lemma yields
$$\liminf_{N\to +\infty}~\P\left(Z'=+\infty,\tau''^N_{\eta}<\sigma''^N_{\eta}\right)\geqslant \P\left(Z'=+\infty\right)-\delta.$$
Therefore (\ref{mino_par_bp}) implies
\begin{align*}\label{milieu_compa}&\liminf_{N\to +\infty}~\P\left(\max_{t\geqslant 0} \Vert I(t)\Vert_1 \geqslant (\alpha-\eta)\Vert z^*\hspace{-1pt}\Vert_1N\right) \notag\\
&\geqslant \liminf_{N\to +\infty}~\P\left(\exists t<\tau''^N_{\eta}\wedge \sigma''^N_{\eta} : \Vert I''(t)\Vert_1 \geqslant (\alpha-\eta)\Vert z^*\hspace{-1pt}\Vert_1N\right)\notag
\\ 
& \geqslant \liminf_{N\to +\infty}~\P\left(Z'=+\infty,\exists t<\tau''^N_{\eta}< \sigma''^N_{\eta} : \Vert I''(t)\Vert_1 \geqslant (\alpha-\eta)\Vert z^*\hspace{-1pt}\Vert_1N\right)\notag
\\
& \geqslant \liminf_{N\to +\infty}~\P\left(Z'=+\infty,\exists t<\tau''^N_{\eta} : \Vert I''(t)\Vert_1 \geqslant (\alpha-\eta)\Vert z^*\hspace{-1pt}\Vert_1N\right)-\delta.\notag
\end{align*}
Proposition \ref{thm_fw_pop} now yields $u>0$ such that 
$$\lim_{N\to +\infty} \P\left(e^{uN}<\tau^N_{\eta}\right)=1$$
so
$$\lim_{N\to +\infty} \P\left(e^{uN}<\tau''^N_{\eta}\right)=1,$$
from which we deduce
\begin{align}
&\liminf_{N\to +\infty}~\P\left(\max_{t\geqslant 0} \Vert I(t)\Vert_1\geqslant (\alpha-\eta)\Vert z^*\hspace{-1pt}\Vert_1N\right)\\ 
& \geqslant \liminf_{N\to +\infty}~\P\left(Z'=+\infty,\exists t<\tau''^N_{\eta} : \Vert I''(t)\Vert_1 \geqslant (\alpha-\eta)\Vert z^*\hspace{-1pt}\Vert_1N\right)-\delta\notag
\\ &
\geqslant \liminf_{N\to +\infty}~\P\left(Z'=+\infty,\exists t<e^{uN} : \Vert I''(t)\Vert_1 \geqslant (\alpha-\eta)\Vert z^*\hspace{-1pt}\Vert_1N\right)-\delta\notag
\\ & \geqslant \P\left((Z'=+\infty)\cap \liminf_{N\to +\infty}~\left(\exists t<e^{uN} : \Vert I''(t)\Vert_1 \geqslant (\alpha-\eta)\Vert z^*\hspace{-1pt}\Vert_1N\right)\right)-\delta\notag
\\ & =\P\left(\liminf_{N\to +\infty}~\left(\exists t<e^{uN} : \Vert I''(t)\Vert_1 \geqslant (\alpha-\eta)\Vert z^*\hspace{-1pt}\Vert_1N\right)\right)-\delta,
\end{align}
where the third inequality proceeds from Fatou's lemma. Now it is well-known (Chapter 1 of \cite{Mode}, Chapter V of \cite{AthNey}, Chapter 4 of \cite{Jag} or \cite{KS67}) that the supercritical branching process $(I''(t))_{t\geqslant 0}$ has a positive exponential growth almost surely on the event where it does not go extinct, so
\begin{equation}\label{fin_compa}\P\left(\liminf_{N\to +\infty}~\left(\exists t<e^{uN} : \Vert I''(t)\Vert_1 \geqslant (\alpha-\eta)\Vert z^*\hspace{-1pt}\Vert_1N\right)\right)=\P\left(\forall t\geqslant 0,I''(t)\neq 0\right).\end{equation}
This last probability is greater than $\P(Z'=+\infty)-\delta$ by definition of $\alpha$, so using (\ref{fin_compa}) then (\ref{milieu_compa}) yields (\ref{a_montrer}), which ends the proof.\hfill $\square$
%
%

\subsection{Proof of Corollary \ref{thm_syst_dyn_exp_dim1}}  

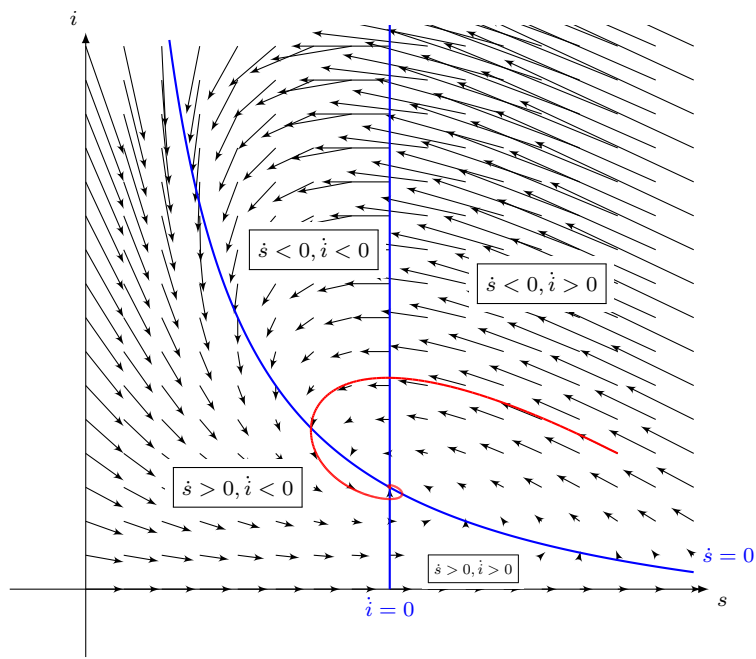
\begin{figure}
\begin{center}
\begin{tikzpicture}[yscale=.9]
\draw [>=latex] [->] (-1,0) -- (8.2,0) node[below right]{$s$};
\draw [>=latex] [->] (0,-1) -- (0,8.2) node[above left]{$i$};
 \clip(-1,-1) rectangle (10,8.3);
	\foreach \u in {0,.5,...,8}
	{\foreach \v in {0,.5,...,8}
		{\draw[thin,>=latex'] [->] (\u,\v) -- +({(10-\u-\u*\v)/20},{(\u*\v-4*\v)/20});};};
\draw[blue, thick] (4,0) node[below, thick, blue]{$\dot i = 0$} -- (4,9);
\draw[thick, domain=1.1:8,samples=200,variable=\x,blue] plot ({\x},{(10-\x)/\x});
\draw[black] (8,.5) node[blue, thick, right]{$\dot s=0$};
\draw (2,1.5) node[fill=white]{\fbox{$\dot s>0, \dot i < 0$}};
\draw (5.1,.3) node[fill=white,scale=.7]{\fbox{$\dot s>0, \dot i > 0$}};
\draw (6,4.5) node[fill=white]{\fbox{$\dot s<0, \dot i > 0$}};
\draw (3,5) node[fill=white]{\fbox{$\dot s<0, \dot i < 0$}};
\edef\debs{7} 
\edef\fins{7} 
\edef\debi{2} 
\edef\fini{2} 
\foreach \t in {1,1.01,...,30} \xdef\debs{\fins} \pgfmathparse{\debs+.01*(10-\debs-\debi*\debs)} \xdef\fins{\pgfmathresult}
\xdef\debi{\fini} \pgfmathparse{\debi+.01*(\debi*\debs-4*\debi)} \xdef\fini{\pgfmathresult} 
\draw[red, thick] (\debs,\debi)--(\fins,\fini)  ;
\end{tikzpicture}
\caption{\label{fig_sir_lv}Projection on the $(s,i)$ plane of the vector field associated with (\ref{simili_lv}) with constant population $z^*=10$, $B=5$, $d=1$, $b=\frac{1}{2}$, $\beta=1$ and $\gamma=3$. Blue curves are nullclines for $s^*$ and $i^*$. The red curve is the solution of (\ref{simili_lv}) starting from point $(7,2)$.}
\end{center}
\end{figure}

It is not difficult to see that $\beta>d+\gamma$ is a necessary and sufficient condition for the existence of an endemic equilibrium for (\ref{simili_lv}) and that the latter is precisely $(s^*,i^*,r^*)$ (see \cite{Nas} for the study of a similar model). If $s(0)+i(0)+r(0)=z^*:=\frac{B}{d-b}$, then the total population is constant and equal to $z^*$ so we can get rid of the third line of $(\ref{simili_lv})$ and the dynamical system can be seen as a Lotka-Volterra prey-predator model (where preys are susceptibles and predators are infectives, see \cite{LV}) with prey immigration. The nullclines for $s$ and $i$ in this model are represented in Figure \ref{fig_sir_lv} along with the associated vector field. For any initial condition in $\R_+\times \R_+^*\times \R_+$, standard arguments (see \cite{LV}) show that $s,i$ and $r$ are well-defined on $\R_+$ and positive. Moreover, $s+i+r$ converges to $s^*+i^*+r^*=z^*$, so it is sufficient to show that $(s,i)$ converges to $(s^*,i^*)$. Setting
$$V(t)=s(t)-s^*\log(s(t))+i(t)-i^*\log(i(t))$$
for all $t\geqslant 0$ yields, after a few calculations:
\begin{align*}\dot{V}(t)&=\underbrace{(B+b(s+i+r)(t))\left(2-\frac{s^*}{s(t)}-\frac{s(t)}{s^*}\right)}_{\leqslant 0}+\underbrace{b\left(z^*-(s+i+r)(t)\right)\left(1-\frac{s(t)}{s^*}\right)}_{=O\left(z^*-(s+i+r)(t)\right)=O\left(e^{(b-d)t}\right)}.&\end{align*}
Now, if $s(0)+i(0)+r(0)\leqslant 2z^*$, $V$ is lower bounded by $-2z^*|\log(4z^*)|>-\infty$, so for any $\eta>0$ its derivative cannot be lower that $-\eta$ for an infinite amount of time. As a result, for any $\delta>0$ we may only have $s(t)\notin [s^*-\delta,s^*+\delta]$ for a finite amount of time because $B+b(s+i+r)>B>0$ and because $\int_0^{+\infty}b\left(z^*-(s+i+r)(t)\right)\left(1-\frac{s(t)}{s^*}\right)\mathrm{d}t$ is finite. This entails that $s$ cannot cross $[s^*-2\delta,s^*-\delta]$ or $[s^*+\delta,s^*+2\delta]$ an infinite number of times since $\dot{s}$ is bounded because of (\ref{simili_lv}). Therefore $s(t)$ lies in $[s^*-\delta,s^*+\delta]$ for $t$ large enough, so $s$ does converge to $s^*$. Similar arguments yield the convergence of $i$ towards $i^*$ using the first equation of (\ref{simili_lv}).
\hfill$\square$

\pe 

%


\bibliography{Biblio}
\bibliographystyle{plain}

\end{document}